\documentclass[a4paper,11pt,english,twoside,leqno]{amsart}

\usepackage[utf8]{inputenc}
\usepackage[T1]{fontenc}
\usepackage[english]{babel}
\usepackage[protrusion=true,expansion=true,kerning, spacing]{microtype}
	\microtypecontext{spacing=nonfrench}
\usepackage{CJKutf8} 
\usepackage{lmodern} 

\usepackage{xspace}

\usepackage[textwidth=16cm,hcentering,heightrounded,foot=1cm]{geometry}
\usepackage{caption}
\usepackage{enumerate,enumitem}
\setlist[enumerate,1]{label=(\roman*)}
\setlist[itemize,1]{leftmargin=20pt}
\usepackage[foot]{amsaddr}

\usepackage{titlesec, titletoc}
\titleformat{name=\section}[block]{\scshape \large \centering}{\thesection.~}{0em}{}[]

\titleformat{\subsection}[hang]{\bfseries}{\thesubsection.~}{0em}{}[]

\titleformat{\subsubsection}[runin]{\bfseries}{\thesubsubsection.}{.2em}{}[]

\titlecontents{section}[0em]{\addvspace{.5pc}\scshape}{\thecontentslabel\  }{}{\titlerule*[.8pc]{.}\contentspage}[]
\titlecontents{subsection}[1.5em]{}{\thecontentslabel\ }{}{\hfill\contentspage}
\titlecontents{subsubsection}[4.5em]{}{\thecontentslabel\ -\ }{}{\hfill\contentspage}

\usepackage{lastpage,fancyhdr}
\usepackage{amsthm,amssymb,amsmath}
\usepackage{amssymb,amsfonts}
\usepackage[fixamsmath]{mathtools} 
\usepackage{accents} 
\usepackage{tikz-cd}
\numberwithin{equation}{subsection}

\usepackage{graphics,graphicx}
\usepackage{color}
\usepackage{tikz}
\usepackage{caption} 
\usepackage{subcaption} 
\usepackage{wrapfig}

\usepackage{csquotes}
\usepackage[
bibencoding=utf8,%
giveninits=true,
bibstyle=ieee-alphabetic,
citestyle=alphabetic,%
sorting=nyt,%
language=english,%
isbn=false,%
doi=false,%
eprint=false,
backref=true, 
hyperref=true,
url=false, 
minbibnames=4,
maxbibnames=4,
dashed=true,
backend=biber,
]{biblatex}
\AtBeginBibliography{\small}

\theoremstyle{definition}
	\newtheorem{defin}{Definition}[section]

\theoremstyle{definition}
	\newtheorem{rem}[defin]{Remark}

\theoremstyle{plain}
	\newtheorem{theo}[defin]{Theorem}
	\newtheorem*{theoU}{Theorem}
	\newtheorem{prop}[defin]{Proposition}
	\newtheorem{cor}[defin]{Corollary}
	\newtheorem{lem}[defin]{Lemma}

\newcommand{\N}{\mathbb{N}}

\newcommand{\Q}{\mathbb{Q}}
\newcommand{\Z}{\mathbb{Z}}

\newcommand{\PP}{\mathbb{P}} 

\newcommand{\B}{\mathbb{B}}

\newcommand{\OO} {{\mathcal O}} 
\newcommand{\mA} { \mathcal{A}}

\newcommand{\mB} {\mathcal{B}}
\newcommand{\mD} {\mathcal{D}}
\newcommand{\mS} {\mathcal{S}}

\newcommand{\mC} {\mathcal{C}}

\newcommand{\mM} {\mathcal{M}}
\newcommand{\mU} {\mathcal{U}}
\newcommand{\mV} {\mathcal{V}}
\newcommand{\mW} {\mathcal{W}}
\newcommand{\mQ} {\mathcal{Q}}
\newcommand{\mfm} {\mathfrak{m}}
\newcommand{\mfX} {\mathfrak{X}}
\newcommand{\mfq} {\mathfrak{q}}
\newcommand{\mfB} {\mathfrak{B}}
\newcommand{\mfA} {\mathfrak{A}}
\newcommand{\mfS} {\mathfrak{S}}

\newcommand{\Ker}{\operatorname{Ker}}
\newcommand{\Spec} { \operatorname{Spec} }
\newcommand{\Spf} { \operatorname{Spf} }

\newcommand{\Hom} {\operatorname{Hom} }

\newcommand{\Card} {\operatorname{Card} }
\newcommand{\Aut} {\operatorname{Aut} }

\newcommand{\Out} {\operatorname{Out} }
\newcommand{\Inn} {\operatorname{Inn} }
\newcommand{\Rev} {\operatorname{Rev} }

\newcommand{\San}{\begin{CJK}{UTF8}{goth}\text{山}\end{CJK}\xspace}
\newcommand{\Ten}{\begin{CJK}{UTF8}{goth}\text{天}\end{CJK}\xspace}
\newcommand{\Zl}{\Z/\ell\Z}
\newcommand{\Zln}{\Z/\ell^n\Z}

\newcommand{\kr}{\underline{kr}}

\newcommand{\Mgm}{\mM_{g,[m]}}
\newcommand{\MgmG}{\Mgm(G)}
\newcommand{\MgmGkr}{\MgmG_{\kr}}
\newcommand{\CompMgmGkr}{\overline{\mM}_{g,[m]}(G)_{\kr}}

\newcommand{\MG}{\Mgm(G)^{\nu}}
\newcommand{\compMG}{\overline{\mM}_{g,[m]}(G)^{\nu}}
\newcommand{\MGkr}{\MG_{\kr}}
\newcommand{\compMGkr}{\compMG_{\kr}}
\newcommand{\MZlkr}{\Mgm(\Zl)^{\nu}_{\kr}}
\newcommand{\compMZl}{\overline{\mM}_{g,[m]}(\Zl)^{\nu}}
\newcommand{\compMZlkr}{\compMZl_{\kr}}

\newcommand{\CgmGkr}{\mathcal{C}_{g,[m]}(G)_{\kr}}
\newcommand{\CgmGGkr}{\mathcal{C}_{g,[m]}[G]_{\kr}}
\usepackage[pdftex,hidelinks,
hyperindex=true, bookmarks=true, bookmarksopen=true, bookmarksnumbered=true,
pdfauthor={Benjamin Collas, Séverin Philip},
pdftitle={On Oda's problem and special loci},
pdfsubject={},
pdfstartview=FitH
]{hyperref}

\begin{document}
	\author{Benjamin Collas, Séverin Philip} 
	\address{Research Institute for Mathematical Sciences, Kyoto University, Kitashirakawa-Oiwakecho, Sakyo-ku, Kyoto 606-8502, Japan}
	\email{bcollas@kurims.kyoto-u.ac.jp, sphilip@kurims.kyoto-u.ac.jp} 
	\urladdr{\url{https://www.kurims.kyoto-u.ac.jp/\textasciitilde bcollas/}}
	\urladdr{\url{https://www.kurims.kyoto-u.ac.jp/\textasciitilde sphilip/}} 
	
	\thanks{Supported by JSPS KAKENHI Grant Number 22F22015 (S.~Philip) and by the Research Institute for Mathematical Sciences, Kyoto University. This work is part of the ``Arithmetic and Homotopic Galois Theory'' project, supported by the CNRS France-Japan AHGT International Research Network between the RIMS Kyoto University, the LPP of Lille University, and the DMA of ENS PSL. Both authors are grateful to Professor A.~Tamagawa for regular discussions on this project and a feedback on a first version of this manuscript.}
	
	\title{On Oda's problem and special loci}
	
	
	\begin{abstract}
		Oda's problem, which deals with the fixed field of the universal monodromy representation of moduli spaces of curves and its independence with respect to the topological data, is a central question of anabelian arithmetic geometry. This paper emphasizes the stack nature of this problem by establishing the independence of monodromy fields with respect to finer special loci data of curves with symmetries, which we show provides a new proof of Oda's prediction.
	\end{abstract}
	
	\maketitle
	\vspace{-3em}
	\tableofcontents
	
	\newpage
	\section{Introduction}
	
	Let $\mM_{g,m}$ be the moduli stack of smooth projective curves of genus $g$ with $m$ (disjoint ordered) sections satisfying the hyperbolicity condition $2g-2+m \geq 1$, which is a smooth geometrically connected Deligne-Mumford stack over $\Q$, and is endowed with a universal punctured curve $\mathcal{C}_{g,m}\to \mM_{g,m}$. For $X$ a punctured curve over $\Q$ of topological type $(g,m)$, associated to a morphism $x\colon\Spec \Q\to \mM_{g,m}$, one obtains two short exact sequences of étale fundamental groups
	\[
	1\to \pi_1^{et}(X\otimes \overline{\Q})\to \pi_1^{et}(\mathcal{C}_{g,m})\to \pi_1^{et}(\mM_{g,m})\to 1\text{ and } 1\to \pi_1^{et}(\mM_{g,m}\otimes\overline{\Q})\to \pi_1^{et}(\mM_{g,m})\overset{p}{\to} G_{\Q}\to 1
	\]
	where the fundamental groups are given with respect to a choice of compatible base points that we omit. Denoting $X\otimes \overline{\Q}$ by $X_{\overline{\Q}}$, the left-hand one gives rise to the \emph{universal $\ell$-monodromy representation}
	\[
	\Phi_{g,m}^{\ell} \colon \pi_1^{et}(\mM_{g,m})\to \Out \pi_1^{et}(X_{\overline{\Q}})\to \Out \pi_1^{\ell}(X_{\overline{\Q}})
	\]
	where the right-hand side morphism comes, for $\ell$ a fixed prime, from the surjective map $\pi_1^{et}(X_{\overline{\Q}}) \to \pi_1^{\ell}(X_{\overline{\Q}})$ to the pro-$\ell$ geometric fundamental group of $X$ (also the maximal pro-$\ell$ quotient of the geometric one). Composing with the section induced by $x$ between Galois and étale fundamental groups, one furthermore recovers the $\ell$-adic representation associated to~$X$
	\[
	\varphi_{X}^{\ell}\colon G_\Q \to \Out \pi_1^{\ell}(X_{\overline{\Q}}).
	\]
	which, contrary to the classical profinite geometric Galois action, has a non-trivial kernel whose corresponding fixed field contains $\Q_{g,m}^{\ell}=\overline{\Q}^{p(\Ker \Phi_{g,m}^{\ell} )}$.
	
	\medskip
	
	The following prediction, as formulated in \cite{iharanak} \S~1.4, stems from Takayuki Oda's original conjecture formulated in \cite{ODAMPI93}. 
	
	\medskip
	
	{\centering
		\begin{minipage}{.9\textwidth}{\bfseries Oda's prediction.} { \itshape For $g$, $m\in\N$ such that $2g-2+m>0$, the $\ell$-monodromy fixed field $\Q_{g,m}^{\ell}$ associated to $\Phi_{g,m}^{\ell}$ is independent of $(g,m)$.}
		\end{minipage}
		\par}
	
	\medskip
	
	As noted in \cite{ODAMPI93}, the group $\Out \pi_1^{\ell}(X_{\overline{\Q}})$ is ``almost intractable'', which motivates Oda to formulate his conjecture in terms of a seemingly more reachable but stronger weight-filtration version of the above prediction, and for fixed $g\geq 0$, see ibid. \emph{\S~2. Theorem and conjectures}. Oda's prediction is finally settled\footnote{Publication of the proof, established in 1995, was indeed postponed to 2012 for unfortunate non-mathematical ground.} for every $(g,m)$ by Takao in \cite{Tak12} following successive advances on the independence in $g$ or $m$ in terms of arithmetic-geometry -- see Ihara and Nakamura in \cite{iharanak}, of group theoretic and Lie algebra computations -- see Nakamura-Takao-Ueno \cite{NTU95} and Matsumoto \cite{matsu} -- and by the use of the (divisorial) Knudsen-Mumford stratification of $\mM_{g,m}$, see \cite{Nacoupl}. An independant proof was later given in terms of combinatorial anabelian geometry by Hoshi and Mochizuki in \cite{HM11}. We also refer to \cite{Tak14} for a recent panorama.
	
	\medskip
	
	Oda's problem -- that is, to which extent canonical arithmetic and geometric data such as $g$ and $m$, produces independent $\ell$-monodromy fixed fields -- is a central question of anabelian arithmetic geometry: it allows the study of the Deligne-Ihara Lie algebra \cite{IHA89} related to motivic multiple zeta values, since for $(g,m)=(0,3)$ the morphism $\Phi_{0,3}^{\ell}$ is the one of Ihara's $\San=\Ten$ question on $\PP_{\Q}^1\setminus\{0,1,\infty\}$ \cite{IHA86}, which in turn, is related to the Rasmussen-Tamagawa conjecture \cite{RASTAM17}. It also has application in low-dimensional topology via the Johnson homomorphism and the Morita obstructions \cite{MOR93}. This conjecture has since motivated the anabelian notion of monodromic fullness \cite{HOS13}.
	
	\medskip
	
	We remark that, as presented in \cite{matsu} Remark~3.3, while Oda's problem is essentially of stack-theoretic nature -- by $\mM_{g,m}$ as a solution to a fine moduli problem and the very existence of the universal punctured curve $\mC_{g,m}$ -- the field $\Q^{\ell}_{g,m}$ was expressed and dealt with in a scheme-theoretic way. This paper develops a setup and techniques that allow to exploit the stack-theoretic aspects of Oda's problem.
	
	\subsection*{Oda's problem for $G$-special loci}
	Let $\Mgm$ denote the moduli stack of curves of genus $g$ with $m$ (unordered) marked points (in particular, $\Mgm$ is not represented by a scheme), which is naturally endowed with a stack inertia stratification, i.e., by the automorphism group of objects. Each strata corresponds to a \emph{$G$-special locus} $\mM_{g,[m]}(G)$ of curves whose automorphism group contains a given finite group $G$. It is shown that the geometric irreducible components $\MgmGkr$ for $G$ cyclic automorphism group, that are among the biggest non-trivial strata, are $\Q$-rational and can be described by combinatorial Hurwitz data $\kr$, see \cite{CM15}.

	\medskip
	
	This context also provides an $\ell$-universal $G$-monodromy representation, see Theorem~\ref{the:univmono}.
	
	\medskip
	
	{\centering
		\begin{minipage}{.9\textwidth}{ \itshape There exists a universal monodromy representation
				\[\Phi^{\ell}_{g,[m]}(G)_{\underline{kr}} \colon \pi_1(\MgmGkr)\longrightarrow \Out \pi_1^{\ell} (X)
				\]
				for $X$ a smooth curve with compactification $\overline{X}$ represented by a $\overline{\Q}$-point on $\MgmGkr$ and $\overline{X}\setminus X$ is a divisor of degree $m$ on $\overline{X}$.
			}
		\end{minipage}
		\par}
	
	\medskip
	
	In particular, this setup provides an $\ell$-monodromy fixed field $\Q_{g,[m]}^{\ell}(G)_{\kr}= \overline{\Q}^{p(\Ker \Phi^{\ell}_{g,[m]}(G)_{\underline{kr}})}$ where $p$ denotes the usual projection to $G_{\Q}$.  In this paper, we deal with the following $\Zln$-special loci version of Oda's problem.
	\medskip
	
	{\centering
		\begin{minipage}{.9\textwidth}{\bfseries Oda's problem for cyclic special loci.} { \itshape For $g$, $m\in\N$ such that $2g-2+m>0$ and $G$ cyclic group of order $\ell^n$, is the $\ell$-monodromy fixed field $\Q_{g,[m]}^{\ell}(G)_{\kr}$ independent of all the special loci data $(g,m)$, $n $ and $\kr$?}
		\end{minipage}
		\par}
	\medskip
	
	While a positive answer to this problem may at first seems ``unreasonable'' -- Oda's problem for cyclic special loci is finer and implies Oda's prediction -- it is supported by a series of indirect results that exhibit \emph{similar arithmetic properties of the stack inertia stratification to the classical divisorial one}: the Galois actions have the same type \cite{CM18}, and the related Grothendieck-Teichmüller groups are isomorphic \cite{Col12}.  More concretely, one notices that the curves used in \cite{matsu} \S~4 to establish Oda's prediction for $2g=0 \mod (\ell-1)$ live in $\MgmGkr$ with $G=\Zl$, quotient genus $g'=0$ and some $\kr$ data with $k=(1,\dots,1,j, -(1+1\cdots +1+j))$ for $j=1$ or $2$, see Section~\ref{subsec:kr} for notations.
	
	\medskip
	
	Indeed, the main results of this paper can be summarized as follows, see Section~\ref{subsec:classical} for the compatibility of the various $\ell$-universal monodromy fields and morphisms and Theorem~\ref{theo:mainFinal}.
	\begin{theoU} Let $\ell$ be a fixed prime. 
		Let $g,m\in \N$ be such that $2g-2+m>0$ and $\underline{kr}$ an associated abstract Hurwitz data such that $\mathcal{M}_{g,[m]}(\Zl)_{\kr}$ is non-empty. 
		The map $\Phi^{\ell}_{g,[m]}(\Zl)_{\underline{kr}}$ is compatible with the map $\Phi^{\ell}_{g,m}$ and the $\ell$-monodromy fixed field $\Q^{\ell}_{g,[m]}(\Zl)_{\underline{kr}}$ is constant equal to $\Q^{\ell}_{0,3}$.
	\end{theoU}

	As a corollary, see Corollary~\ref{cor:newOda}, we recover the containment $\Q^{\ell}_{g,m}\subset \Q^{\ell}_{0,3}$ and thus the classical version of Oda's prediction, that is for all $g$, $m\in \N$ such that $2g-2+m>0$ we have $\Q^{\ell}_{g,m}= \Q^{\ell}_{0,3}$, see also \cite{iharanak} Theorem~3~B. Both proofs of Oda's problem for special loci and classical settings still rely on the previously established $\Q^{\ell}_{0,3}\subset \Q^{\ell}_{g,m}$, see \cite{Nacoupl,matsu,Tak12}.
	
	\medskip
	
	\emph{The organization of the paper is as follows.} In Section~\ref{sec:OdaConjG} we recall the $\kr$ combinatorial description of irreducible components of cyclic special loci of \cite{CM15} and introduce the $\ell$-universal $G$-monodromy representation, whose fixed field we relate within a lattice of other $\ell$-monodromy fixed fields, which in particular includes the more traceable Hurwitz spaces $\Mgm[G]_{\kr}$
	\[
	\begin{tikzcd}[column sep=10pt]
		\Q^{\ell}_{g,m} \arrow[r, hook] &  \Q^{\ell}_{g,[m]}(\Zln)_{\underline{kr}} \arrow[r, hook] & \Q^{\ell}_{g,[m]}(\Zln)_{\underline{kr}}^\nu\\
		\Q^{\ell}_{0,3} \arrow[u, hook, dashed] \arrow[r, hook,dashed] & \Q^{\ell}_{g',m'}  \arrow[r, hook]&  \Q^{\ell}_{g',[m']}(\delta \Zln)^\nu \arrow[u, hook] \ar[u, shorten=2mm, dash, shift left]
		&
	\end{tikzcd}
	\]
	where $(g',m')$, resp.$\Q^{\ell}_{g',[m']}(\delta \Zln)^\nu$, denotes the topological data, resp. a certain monodromy fixed field, obtained by $G$-quotient. At this stage, establishing the $G$-special version of Oda's prediction relies on showing that $\Q_{g,m}^{\ell}(\Zln)_{\kr}^\nu\subset \Q_{0,3}^{\ell}$; our proof adapts Ihara-Nakamura's \cite{iharanak}. Section~3 deals with the construction of tangential base points, or one-parameter families, on the $G$-stable compactification of Hurwitz spaces in terms of formal patching of certain Matsumoto-Seyama $\Zl$-stable curves, whose Galois action properties are established in Section 4 via Grothendieck-Murre theory and by comparison with Deligne's original tangential base point. This results in the inclusion of the $\ell$-monodromy fixed field of the generic fiber of the constructed one-parameter families into $\Q_{0,3}^{\ell}$. We conclude with a general Theorem~\ref{theo:GM4App} that can be applied to multiple geometric situations. Section~5 ties everything together for $\Zl$, starting with the case of proper loci for which the deformation method does not apply. In the diagram above, Oda's classical prediction then follows the bottom row reading.

	\bigskip
	
	\emph{Notations and conventions.} For $G$ a finite group, we write $\mM_{g,[m]}[G]$ for the Hurwitz space of $G$-covers and $\MgmG^{\nu}$ for the quotient $\mM_{g,[m]}[G]/\Aut G$. We denote by $\overline{\mM}_{g,[m]}(G)$ the $G$-stable compactification of the $G$-special locus $\MgmG$, and by $\compMG$ the stable compactification of $\MgmG^{\nu}$.  The topological data $(g,m)$ of a curve are said to be of hyperbolic type if they satisfy $2g-2+m>0$.

		\section{Oda's conjecture for $G$-special loci}\label{sec:OdaConjG} 
		After some brief reminders on the description of irreducible components $\MgmGkr $ of cyclic $G$-special loci in terms of combinatorial Hurwitz data $\kr$, we define the $\ell$-universal $G$-monodromy representation $\Phi^{\ell}_{g,m}(G)_{\underline{kr}}\colon \pi_1(\MgmGkr) \longrightarrow \Out \pi_1^{\ell}(X)$ -- for $G$ any finite group -- where $X$ is a hyperbolic curve of type $(g,m)$. Relying on the forgetful functor and the quotient functor
		\[
		\MgmGkr^{\nu} \rightarrow \MgmGkr \rightarrow \Mgm, \text{ and } \mM_{g,[m]}[G]_{\underline{kr}}  \overset{\delta}{\rightarrow} \mM_{g',[m']}	
		\]
		and some properties of the stack inertia $\mathcal{I}_\mM$, we build step-by-step \emph{a lattice of relations between the various $\ell$-monodromy fixed fields arising from this context} -- that is between $\Q_{g,m}^{\ell}$, $\Q_{g,[m]}^{\ell}(\Zln)_{\kr}$, $\Q_{g,[m]}^{\ell}(\Zln)_{\kr}^\nu$, and $\Q_{g',[m']}^{\ell}(\delta\Zln)_{\kr}^\nu$. 
		
		\subsection{Universal monodromy representations and Oda's fields for $G$-special loci}\label{sec:MonOdaG}
		\subsubsection{}\label{subsec:kr}
		Let $\Mgm[G]$ be the moduli stack of curves of genus $g$ with $m$ marked points endowed with a faithful $G$-action, or \emph{Hurwitz stack}, whose $S$-sections for a $\Q$-scheme $S$ are defined as follows: 
		\[
		\Mgm[G](S) \text{ are the triplets } (C,D,\iota)
		\text{ where }
		\begin{cases}
			C\text{ is a smooth projective curve of genus $g$ over }S,\\
			\iota \colon G\rightarrow \Aut_S C\text{ an injective homomorphism,}\\
			D\text{ an étale Cartier divisor of degree } m \\
			\quad \text{ stabilized by the } G\text{-action induced by } \iota,
		\end{cases}
		\]
		see \cite{CM15} \S~2.1 as well as for the rest of this section. The \emph{$G$-special locus $\Mgm(G)$ of $\Mgm$} is obtained as the image of $\Mgm[G]$ in $\Mgm$ under the forgetful functor defined by
		\[
		\begin{array}{cccc}
			\Mgm[G](S) & \longrightarrow & \Mgm(S)  \\
			(C,D,\iota) & \longmapsto & (C,D)
		\end{array}
		\]
		In particular, the $S$-sections of $\Mgm(G)$ are curves over $S$ whose geometric fibers admits a faithful $G$-action. The stack $\Mgm[G]$ having a canonical right-action of $\Aut G$ via $\iota$, we can form the quotient stack $\Mgm[G]/\Aut G$ that we denote by $\MG$ since, apart from a few exceptional cases\footnote{Erratum: Proposition~2.4 and Corollaire~2.5 of \cite{CM15} are subject to the same exceptions.} see \cite{MSSV02} Theorem~5.1 and section~4 for an account with $g\geq 2$ and also Remark~\ref{rem:HurNonNorm} \ref{it:nonNormLoc}, it identifies with the normalization of $\Mgm(G)$ by the proof of \cite{Ro11} Proposition~3.4.1. All the stacks $\Mgm[G]$, $\Mgm(G)$ and $\MG$ are Deligne-Mumford stacks over $\Spec \Q$ -- with $\Mgm[G]$ and $\MG$ moreover smooth over $\Spec \Q$.
		
		\medskip
		
		From now on, \emph{we assume that $G\simeq \Z/n\Z$ is cyclic}, so that following \cite{CM15}, we can investigate the subloci $\MgmGkr$ of $\Mgm(G)$ of $S$-curves whose $G$-action ramification data correspond to certain Hurwitz data $\underline{kr}=(k,r)$ modulo the diagonal $(\Z/n\Z)^{\times}$-action, which are abstractly defined as follows:
		\begin{itemize}
			\item The part $k$ corresponds to an $N$-tuple in $(\Z/n\Z)^{N}$, where $N$ is the degree of the branch divisor, whose terms sum to $0$, and which is taken up to permutation. Each component of $k$ corresponds to a generator of one of the $G$-isotropy groups.
			\item The second part $r$ is an element of $\N^n$, whose $i$-th component, in the case of a quotient map $\psi\colon C\rightarrow C/G$, corresponds to
			\[ 
			r(i)= \Card \{ y\in D/G \mid br(y)= i \mod n\} 
			\]
			where $br(y)$ is the branching order at $y$, that is the ramification index of any point in the fiber $\psi^{-1}(y)$.
		\end{itemize}
		Note that the $(\Z/n\Z)^{\times}$-quotient in $\underline {kr}$ should be seen as the $(\Aut G)$-quotient previously introduced. We refer to ibid. Definitions~3.5 and~3.9, and Example~3.11 for further details.
		
		\medskip
		
		The construction of abstract Hurwitz data from $G$-curves defines a map
		\[
		\underline{kr}\colon  \Mgm[G]_{N}  \longrightarrow ((\Z/n\Z)^{N}/\mathfrak{S}_{N}\times \N^n)/(\Z/n\Z)^{\times}   
		\]
		where $\Mgm[G]_{N}$ denotes the substack of $\Mgm[G]$ of curves whose branch divisor is of degree $N$, which is locally constant -- see \cite{CM15} Lemma~3.13. For a fixed value of $\underline{kr}$, one thus obtains a substack $\Mgm[G]_{\underline{kr}}$ of $\Mgm[G]$ of $G$-curves with abstract Hurwitz data $\underline{kr}$ so that one can define:
		
		\begin{defin} For $G$ cyclic and given abstract Hurwitz data $\underline{kr}$ the special sublocus $\MgmGkr$ is the image of $\Mgm[G]_{\underline{kr}}$ under the forgetful functor $\Mgm[G]\rightarrow \Mgm(G)$.
		\end{defin}
		
		Also, since the action of $\Aut G$ stabilizes $\Mgm[G]_{\underline{kr}}$ by definition, we have substacks $\MgmGkr^{\nu}$ of $\MgmG^{\nu}$. The stacks $\MgmGkr^{\nu}$ and $\MgmGkr$ are defined over $\Q$ by construction and are geometrically irreducible by Proposition~3.12 and Theorem~4.3 of \cite{CM15}.
		
		\medskip
		
		One particular case of interest is when the ramification divisor is contained in the marked divisor $D$. In this case, we can recover $r$ by the data of $D$ and $k$. Indeed, we have 
		\[
		\begin{cases}
			r(i)=  \Card \{j \mid k(j)=i\}/ \operatorname{gcd}(i,n) \text{ for } i\neq 0\\ 
			r(0)= \deg D - \sum_{i\in \Z/n\Z\setminus \{0\}} \Card \{j \mid k(j)=i\}.
		\end{cases}
		\]
		
		\medskip
		
		Similarly to the moduli stacks of curves, the stacks $\mM_{g,[m]}(G)$, resp. $\MG $, are not necessarily proper. We denote by $\overline{\mM}_{g,[m]}(G)_{\kr}$ the $G$-stable compactification of the $G$-special locus $\MgmGkr$, and by $\compMG_{\kr}$ the $G$-stable compactification of $\MGkr$. These are obtained from the original stacks by adding stable curves endowed with a stable $G$-action. We refer to \cite{EKE95} and \cite{BERO07} \S~4 and~6 for details.

		\begin{rem}\label{rem:HurNonNorm}\mbox{}
			\begin{enumerate}
				\item The correspondence between the abstract Hurwitz data $\kr$ and the Hurwitz data $\xi$ of \cite{BERO07} \S~2.2 in terms of equivalence classes $[H_i,\chi_i]$ of characters $\chi_i$ at $G$-inertia group $H_i$ is straightforward by considering generators of the $G$-isotropy groups.
				\item\label{it:nonNormLoc} The difference between $\mM_{g,[m]}(G)$ and $\MG $ comes from the potential existence of a curve whose geometric fiber has an automorphism group that contains 2 topologically but not holomorphically conjugate subgroups. We refer to \cite{GDHIL97} for examples.
			\end{enumerate}
		\end{rem}

		\subsubsection{}
		We now consider $\CgmGkr$ the universal $G$-curve of genus $g$ with $m$ punctures and abstract Hurwitz data $\underline{kr}$. We denote by $\mM_{g,[m]+1}$ the stack of smooth projective curves with a degree $m$ divisor and an additional marked point. We have an identification $\CgmGkr\simeq \MgmGkr\times_{\Mgm} \mM_{g,[m]+1}$. The $S$-sections of $\CgmGkr$ are the elements of $\MgmGkr(S)$ with the additional data of a section outside the marked points $D$; similarly, the universal punctured curve over $\MgmGkr^{\nu}$ is given by the stack $\CgmGkr^{\nu}=\MgmGkr^{\nu}\times_{\Mgm} \mM_{g,[m]+1}$.
		
		\medskip
		
		One obtains the $\ell$-universal $G$-monodromy representation. 
		\begin{theo} \label{the:univmono} Let $g$, $m \in \N$ such that $2g-2+m>0$, $G$ a finite cyclic group and $\kr$ a Hurwitz data with respect to $g$, $m$ and $G$, then there is an exact sequence
			\[\begin{tikzcd}
				1 \arrow[r] &  \widehat{F}_{2g+m-1} \ar[r] & \pi_1(\CgmGkr) \arrow[r] &  \pi_1(\MgmGkr) \arrow[r] & 1.
			\end{tikzcd}\]
			The \emph{$\ell$-universal $G$-monodromy representation} is the induced monodromy map 
			\begin{equation}\label{eq:UniMonoGloc}
				\Phi^{\ell}_{g,m}(G)_{\underline{kr}}\colon \pi_1(\MgmGkr) \longrightarrow \Out \widehat{F}^{\ell}_{2g+m-1} 
			\end{equation}
			which is universal in the following sense: for any curve $C$ over a connected $\Q$-scheme $S$ in $\MgmGkr(S)$ and $\overline{\Q}$-point $\overline{s}$ of $S$ the natural representation $\pi_1(S)\rightarrow \Out \pi^{\ell}_1(C_{\overline{s}})$ factors through $\Phi^{\ell}_{g,m}(G)_{\underline{kr}}$. A similar result holds for $\MgmGkr^{\nu}$.
			
		\end{theo}
		
		In the exact sequence above one has identified the fundamental group of the fiber of the map $\CgmGkr\rightarrow \MgmGkr$ at the geometric base point with $\widehat{F}_{2g+m-1}$. In the same way, the factorization of the representation to $\Out \pi^{\ell}_1(C_{\overline{s}})$ through $\Phi^{\ell}_{g,m}(G)_{\underline{kr}}$ is made via the identification $\pi^{\ell}_1(C_{\overline{s}})\simeq\widehat{F}^{\ell}_{2g+m-1}$.

		\begin{proof}  Let $\overline{x}\colon \Spec \overline{\Q} \rightarrow \MgmGkr$ be a geometric point representing a curve $X$ over $\overline{\Q}$. By taking the rigidification given by a Jacobi structure of level $N\geq 3$ we obtain étale Galois covers $\MgmGkr^N$ and $\CgmGkr^N$ of $\MgmGkr$ and $\CgmGGkr$, respectively, which are schemes and sit in a similar sequence, compare with \cite{DM69} \S~5.4 and \S~5.14. The induced maps from this new sequence to the old one make the following commutative diagram, with exact columns and bottom row,
			\[\begin{tikzcd}
				1 \arrow[r] &  \pi_1(X) \ar[r] \arrow[d, equal] & \pi_1(\CgmGkr^N) \arrow[r] \arrow[d, hook] &  \pi_1(\MgmGkr^N) \arrow[r] \arrow[d, hook]  & 1 \\
				1 \arrow[r] &  \pi_1(X) \ar[r] \arrow[d] & \pi_1(\CgmGkr) \arrow[r] \arrow[d, twoheadrightarrow]&  \pi_1(\MgmGkr) \arrow[r] \arrow[d, twoheadrightarrow] & 1 \\
				& 1 \arrow[r] &  \Aut (\Z/N\Z)^{2g}  \arrow[r] & \Aut (\Z/N\Z)^{2g} \arrow[r]  &  1
			\end{tikzcd}\]
			
			By a diagram chase the exactness of the upper sequence implies that of the middle one. The right exactness of the upper sequence is given by \cite{SGA1} Exposé~IX Corollaire~6.11. The left exactness then follows from the hyperbolicity condition and the identification with the profinite completion of the Birman exact sequence.

			A similar argument provides the result for $\MGkr$ with \emph{ad hoc} substitutions.
		\end{proof}

		For a curve $C$ over $S$ as in Theorem~\ref{the:univmono} the \emph{$\ell$-monodromy representation of $C$}
		\[
		\varphi^{\ell}_C\colon \pi_1(S)\rightarrow \Out \pi_1^{\ell} (C_{\overline{s}})
		\]
		is obtained from the relative homotopy exact sequence as usual. Notice that the $\Q$-scheme $S$ also sits in a classical arithmetic-geometric homotopy exact sequence, so that $\pi_1(S)$ is naturally equipped with a projection map $p_S\colon \pi_1(S)\rightarrow G_{\Q}$. We recall that, similarly, we have a canonical homomorphism $p\colon \pi_1(\MgmGkr)\to G_{\Q}$.
		
		\begin{defin}
			The field $\Q^{\ell}_{g,[m]}(G)_{\underline{kr}}$, resp. $\Q^{\ell}_{g,[m]}(G)^{\nu}_{\underline{kr}}$, is the fixed field of $p(\Ker\Phi^{\ell}_{g,[m]}(G)_{\underline{kr}})$, resp. of $p(\Ker\Phi^{\ell}_{g,[m]}(G)_{\underline{kr}}^\nu)$. 
			For a curve $C$ over a connected $\Q$-scheme $S$, the field $\Q^{\ell}_C$ is the fixed field of $p_S(\Ker \varphi^{\ell}_C)$.
		\end{defin}

		\begin{lem}\label{lem:inclusioncurve}
			For $C$ a curve over a connected $\Q$-scheme $S$ represented by an $S$-point on $\MgmGkr$, resp. on $\MgmGkr^{\nu}$, we have the inclusion
			\[
			\Q^{\ell}_{g,[m]}(G)_{\underline{kr}}\subset \Q^{\ell}_C,\text{ resp. } \Q^{\ell}_{g,[m]}(G)_{\underline{kr}}^\nu\subset \Q^{\ell}_C.
			\]
			The $\ell$-monodromy fixed field $\Q^{\ell}_{g,[m]}(G)_{\underline{kr}}$ is furthermore obtained as the intersection of all the $\Q^{\ell}_C$ for such $C/S$ where $S$ varies in the category of connected $\Q$-schemes.
		\end{lem}
		
		The field $\Q^{\ell}_{g,[m]}(G)_{\underline{kr}}$ can also be obtained as $\Q^{\ell}_{C_0}$ where $C_0= \CgmGkr \times_{\MgmGkr}\mM_{g,[m]+m'}(G)_{\underline{kr}}$ is a curve over $S=\mM_{g,[m]+m'}(G)_{\underline{kr}}$ with $m'$ large enough for $S$ to be a scheme.
		
		\begin{proof}
			By the universality of the map $\Phi^{\ell}_{g,[m]}(G)_{\underline{kr}}$ we have a commutative diagram
			\[\begin{tikzcd}
				\pi_1(S) \arrow[d, "p_S"'] \arrow[r]  & \pi_1(\MgmGkr) \arrow[r] \arrow[d, "p"]& \Out \widehat{F}^{\ell}_{2g+m-1}\\
				G_{\Q} \ar[r, equal] & G_{\Q} & \\
			\end{tikzcd}
			\]
			where $\varphi^{\ell}_C$ appears as the composition $\pi_1(S)\rightarrow \pi_1(\MgmGkr) \rightarrow \Out \widehat{F}^{\ell}_{2g+m-1}$. The compatibility with the projections to $G_{\Q}$ ensures that we have $p_S(\Ker \varphi^{\ell}_C) \subset p(\Ker \Phi^{\ell}_{g,[m]}(G)_{\underline{kr}})$ and thus the inclusion. To prove the last point, by commutativity of the diagram, it suffices to show the existence of a curve $C$ in $\MgmGkr(S)$ such that the induced map $\pi_1(S)\rightarrow \pi_1(\MgmGkr)$ is surjective. This is done by taking $C_0= \CgmGkr \times_{\MgmGkr}\mM_{g,[m]+m'}(G)_{\underline{kr}}$ over $S=\mM_{g,[m]+m'}(G)_{\underline{kr}}$ with $m'$ large enough for $S$ to be a scheme.
			
			\medskip
			
			The case of $\Q^{\ell}_{g,[m]}(G)_{\underline{kr}}^\nu$ is similar after replacing $\MgmGkr$ by $\MGkr$.
		\end{proof}

		\subsubsection{}\label{sub:kret}
		Let us now relate the general situation to the one where the divisor of marked points $D$ contains the ramification divisor $R$ of the $G$-action, a property that we recall, can be seen directly on the abstract Hurwitz data.
		
		\medskip
		By base change to an algebraically closed field and reading of the $\underline{kr}$ data one notices that the divisor $R\cup D$ is finite étale over $S$ for a curve $C/S$ as before.
		
		\begin{lem}
			Let $(C,D)$ be a curve represented by an $S$-point on $\MGkr$ as before. Then the degree of the ramification divisor $R$ of $C$ and of the divisor $R\cup D$ are determined by the abstract Hurwitz data $\underline{kr}$.  
		\end{lem}
		
		\begin{proof}
			As everything is locally constant on the base, it is enough to treat the case where $S$ is the spectrum of an algebraically closed field. By definition of $\underline{kr}$ the degree $\deg R=N$ of the ramification divisor is the length of $k$. Furthermore, since the degree of $R\cap D$ is given by $\sum\nolimits_{i=1}^{n-1} \operatorname{gcd}(i,n)\cdot \underline{r}(i)$, we have the formula
			\[
			\deg R\cup D=m+N-\sum\limits_{i=1}^{n-1} \operatorname{gcd}(i,n)\cdot r(i)
			\]
			which is entirely determined by $m$, $\underline{kr}$ and $G=\Z/n\Z$.
		\end{proof}
		
		For an abstract Hurwitz data $\underline{kr}$, we introduce $\underline{kr}^{et}$ as the minimal associated Hurwitz data such that the ramified points are contained in the marked divisor -- i.e. minimal in the sense that the new marked divisor is the smallest one containing $D$ and $R$ -- and which is thus defined by
		\[\begin{cases}
			r^{et}(0)=r(0) &\\
			r^{et}(i)=\Card \{j\in \{1,\dots, N\} \mid k(j)=i\}, & i\geq 1.
		\end{cases}\]

		\begin{prop}\label{pro:etmap}
			There is a natural map of stacks 
			\[
			\MGkr\longrightarrow \mM_{g,[m+s]}(G)_{\underline{kr}^{et}}^{\nu}
			\]
			where $\underline{r}^{et}$ and $s=\deg R-\deg R\cap D$ can be explicitly determined as above.
		\end{prop}
		
		\begin{proof}
			By the previous lemma we have that if $(C,D)$ is in $\mM_{g,[m]}(G)_{\underline{kr}}^{\nu}(S)$ then $(C,R\cup D)$ is an element of $\mM_{g,[m+s]}(G)_{\underline{kr}^{et}}^{\nu}(S)$. This association defines a map of groupoids as any isomorphism preserving the $G$-action must also preserve the ramification divisor. 
		\end{proof}

		\begin{theo}\label{the:nonetale}
			We have the following inclusion of $\ell$-monodromy fixed fields
			\begin{equation}
				\Q^{\ell}_{g,[m]}(G)^\nu_{\underline{kr}}\subset \Q^{\ell}_{g,[m+s]}(G)^\nu_{\underline{kr}^{et}}.
			\end{equation}
		\end{theo}
		
		\begin{proof}
			Let $\sigma \in p(\Ker \Phi^{\ell}_{g,[m+s]}(G)_{\underline{kr}^{et}})\subset G_{\Q}$. By Lemma~\ref{lem:inclusioncurve} there is a connected $\Q$-scheme $S$ and a curve $(C,D)$ over $S$ represented by an $S$-point on $\mM_{g,[m+s]}(G)_{\underline{kr}^{et}}$ such that $\sigma$ has a lift $\tau$ in the kernel of the map
			\[
			\begin{tikzcd}
				\pi_1(S)\arrow[r, "s_C"] & \pi_1(\mM_{g,[m+s]}(G)_{\underline{kr}^{et}}) \arrow[r, "\Phi^{\ell}_{g,[m+s]}(G)_{\underline{kr}^{et}}"] &[30pt] \Out(\widehat{F}_{2g+m+s-1}^{\ell}). 
			\end{tikzcd}
			\]
			The divisor $D$ admits a decomposition $D=D^{un}\cup D^{ram}$ where $D^{un}$ is given by the unramified marked points and $D^{ram}$ by the ramified marked points. By definition of the component $\underline{r}^{et}$ of $\kr^{et}$, the divisor $D^{ram}$ corresponds to all the ramified points. The divisor $D^{ram}$ splits into a disjoint union of geometrically irreducible divisors over a finite étale extension $S'=S_K$ of $S$ where $K$ is defined by the property that $G_K$ stabilizes each geometric component of $D^{ram}$. In particular, $\pi_1(S')$ contains the subgroup $\{ \alpha \in \pi_1(S) \mid p_S(\alpha)\in G_K\}$, which contains $\tau$ by construction. We can thus assume that $S=S'$.
			
			By removing some chosen orbits of ramified points in $D^{ram}$ according to the data given by $\underline{r}$ we can form a divisor $D'=D^{un}\cup D^{ram'}$ such that $(C,D')$ gives an $S$-point of $\mM_{g,[m]}(G)_{\underline{kr}}$. Hence, it is sufficient to show that $\sigma$ is the image of an element of $\pi_1(S)$ that acts trivially on the pro-$\ell$-fundamental group of a geometric fiber $C_{\overline{s}}\setminus D'_{\overline{s}}$ of $C\setminus D'$. This now comes from the fact that the outer actions of $\pi_1(S)$ on $\pi_1^{\ell} (C_{\overline{s}}\setminus D_{\overline{s}})$ and $\pi_1^{\ell}(C_{\overline{s}}\setminus D'_{\overline{s}})$ are compatible with the canonical surjection $\pi_1^{\ell} (C_{\overline{s}}\setminus D_{\overline{s}})\rightarrow \pi_1^{\ell}(C_{\overline{s}}\setminus D'_{\overline{s}})$.  
		\end{proof}

		\subsection{From the classical to the special loci settings} \label{subsec:classical}
		In order to relate the $\ell$-monodromy fixed fields $\Q_{g,m}^{\ell}$ and $\Q_{g,m}^{\ell}(G)_{\kr}$ let us start by showing that we can move from $\mM_{g,m}$ to $\mM_{g,[m]}$ without harm. Let $\Q^{\ell}_{g,[m]}$ be the fixed field of $p(\Ker \Phi^{\ell}_{g,[m]})$ where $p\colon \mM_{g,[m]}\rightarrow \Spec\Q$ is the structure map and $\Phi^{\ell}_{g,[m]}\colon \pi_1(\mM_{g,[m]}) \rightarrow \Out \pi_1^{\ell}(C)$ the outer Galois action coming from the exact sequence
		\[\begin{tikzcd}
			1 \arrow[r] & \pi_1(C) \arrow[r] & \pi_1(\mM_{g,[m]+1}) \arrow[r] & \pi_1(\mM_{g,[m]}) \arrow[r] & 1
		\end{tikzcd}
		\]
		where $C$ is a geometric fiber of $\mM_{g,[m]+1}\rightarrow \mM_{g,[m]}$ . The following can also be seen as a special case of \cite{Ho11} Lemma 1.4~(ii). 
		\begin{lem}\label{lem:nobrack}
			We have $\Q^{\ell}_{g,m}=\Q^{\ell}_{g,[m]}$. 
		\end{lem}
		
		\begin{proof}
			It suffices to see that the equality $\Ker \Phi^{\ell}_{g,m}=\Ker \Phi^{\ell}_{g,[m]}$ holds in $\pi_1(\Mgm)$ as we have $\pi_1(\mM_{g,m})\subset \pi_1(\Mgm)$ with cokernel $\mathfrak{S}_m$. For a presentation of $\pi_1^{\ell}(C)$ given by
			\[
			\langle y_1,\dots,y_{2g}, x_1,\dots, x_m \mid [y_1,y_2]\cdots [y_{2g-1},y_{2g}]x_1\cdots x_m=1\rangle
			\]
			it is clear that an element $\tau\in \pi_1(\Mgm)$ has image $\sigma\in \mathfrak{S}_m$ if and only if the permutation induced by $\tau$ on the set of conjugacy classes of cuspidal inertia subgroups of $\pi_1^{\ell}(C_{\overline{\Q}})$, which is in bijection with the set $\{x_1,\dots, x_m\}$,  is the one given by $\sigma$. Such an element $\tau$ thus has trivial outer action on $\pi_1^{\ell}(C)$ only if it has trivial image in $\mathfrak{S}_m$ and thus belongs to $\pi_1(\mM_{g,m})$.
		\end{proof}
		
		\subsubsection{ }The comparison via the forgetful functor $\MgmGkr^{\nu} \rightarrow \MgmGkr \rightarrow \Mgm$ is now straightforward.
		
		\begin{prop}\label{pro:compgn}
			For all $(g,m)$ of hyperbolic type and compatible Hurwitz data $\kr$ we have $\Q_{g,m}^{\ell} \subset \Q_{g,[m]}^{\ell}(G)_{\kr}\subset  \Q_{g,[m]}^{\ell}(G)_{\kr}^\nu$.
		\end{prop}
		
		\begin{proof}
			Let $C$ be a curve over $\overline{\Q}$ represented on $\MGkr$. 
			First see that the sequence of maps
			\[
			\pi_1(\MGkr) \rightarrow \pi_1(\MgmGkr) \rightarrow \pi_1(\Mgm) \rightarrow \Out \pi_1^{\ell}(C)
			\]
			induces a sequence
			\[
			\Ker \Phi^{\ell}_{g,[m]}(G)_{\underline{kr}}^\nu \rightarrow \Ker \Phi^{\ell}_{g,[m]}(G)_{\underline{kr}} \rightarrow  \Ker \Phi^{\ell}_{g,[m]}
			\]
			where the second map is obtained by considering the following commutative diagram with exact rows  
			\[ \begin{tikzcd}
				1 \arrow[r] & \pi_1(C) \arrow[d, equal] \arrow[r] & \pi_1( \CgmGkr) \arrow[r] \arrow[d] & \pi_1( \MgmGkr) \arrow[d] \arrow[r] & 1 \\
				1 \arrow[r] & \pi_1(C) \arrow[r]\arrow[d] & \pi_1(\mM_{g,[m]+1}) \arrow[r]\arrow[d] & \pi_1(\mM_{g,[m]}) \arrow[d] \arrow[r] & 1\\
				1 \arrow[r] & \Inn \pi_1^{\ell}(C) \arrow[r] & \Aut \pi_1^{\ell}(C) \arrow[r] & \Out \pi_1^{\ell} (C)  \arrow[r] & 1
			\end{tikzcd} \]
			and the first map is obtained in a similar way.
			
			By applying the canonical projections to $G_{\Q}$, and Lemma~\ref{lem:nobrack} for $\Q^{\ell}_{g,m}=\Q^{\ell}_{g,[m]}$, one obtains the desired sequence of inclusions.
		\end{proof}
		
		\begin{cor}\label{cor:nonetale}
			With the notations of Theorem~\ref{the:nonetale} we have 
			\[
			\Q^{\ell}_{0,3} \subset \Q_{g,[m]}^{\ell}(G)_{\kr}\subset  \Q_{g,[m]}^{\ell}(G)_{\kr}^\nu \subset \Q_{g,[m+s]}^{\ell}(G)_{{\kr}^{et}}^\nu.
			\]
		\end{cor}
		\begin{proof}
			The inclusion $\Q_{0,3}^{\ell}\subset \Q^{\ell}_{g,m}$ for all hyperbolic $(g,m)$ is essentially Theorem~3.6 of \cite{Tak12}. The rest of the inclusions follow from Proposition~\ref{pro:compgn} and Theorem~\ref{the:nonetale}. 
		\end{proof}
		
		\begin{rem}\label{rem:OdaWeightOK}
			In Proposition~\ref{pro:compgn} there is no difficulty to move to the weight version of Oda's conjecture, and we get, for all $(g,m)$ of hyperbolic type, any compatible Hurwitz data $\underline{kr}$, and all weight $w$
			\[
			\Q_{g,m}^{\ell}(w) \subset \Q_{g,[m]}^{\ell}(G)_{\underline{kr}}(w)\subset  \Q_{g,[m]}^{\ell}(G)_{\underline{kr}}^\nu(w).
			\] 
			In contrast see Remark~\ref{rem:beyondOdanoWeight}~\ref{it:rem:OdaWeightProof}.	
		\end{rem}
		
		\subsubsection{} The quotient map $\delta\colon \mM_{g,[m]}[G]_{\underline{kr}} \to \mM_{g',[m']}$ defined by $(C,D,\iota) \mapsto (C/\iota(G),D/\iota(G))$ allows the comparison of $\ell$-monodromy fixed fields. We first remark that the map $\delta$ is well-defined at the level of the stack $\MGkr$, since $\delta$ is equivariant under the action of $\Aut G$.
		
		\medskip

		Therefore we have a map $\delta\colon \MGkr\rightarrow \mM_{g',[m']}$ that fits in a commutative square
		\[
		\begin{tikzcd}
			\CgmGkr^{\nu}  \arrow[r] \arrow[d] & \MgmGkr^{\nu} \arrow[d, "\delta"] \\
			\mM_{g',[m']+1} \arrow[r] & \mM_{g',[m']}
		\end{tikzcd}
		\]  
		where the map on the left is induced by the quotient in the same way. For a curve $\overline{X}$  over $\overline{\Q}$ represented on $\MgmGkr^{\nu}$ let us denote $\overline{Y}$ the quotient proper curve, and $X$, $Y$ their open counterparts. This leads to a commutative diagram with exact rows
		
		\[\begin{tikzcd}
			1\arrow[r] & \pi_1(X) \arrow[r] \arrow[d]& \pi_1(\CgmGkr^{\nu}) \arrow[r] \arrow[d] & \pi_1(\MgmGkr^{\nu}) \arrow[d] \arrow[r] & 1  \\
			1\arrow[r] & \pi_1(Y) \arrow[r] & \pi_1(\mM_{g',[m']+1}) \arrow[r] & \pi_1(\mM_{g',[m']}) \arrow[r] & 1 \\
		\end{tikzcd}
		\]
		which in turn provides an $\ell$-monodromy representation 
		\[
		\Phi_{g',[m']}^{\ell}(\delta G)_{\underline{kr}}^\nu\colon\pi_1(\MgmGkr^{\nu}) \rightarrow \Out \pi_1^{\ell}(Y)
		\] 
		in the quotient curve, so that one obtains
		\begin{equation}\label{Eq:gmDetlaFields}
			p(\Ker \Phi_{g',[m']}^{\ell}(\delta G)^\nu_{\underline{kr}})\subset p(\Ker \Phi_{g',[m']}^{\ell})\text{ or equivalently } \Q^{\ell}_{g',[m']}\subset \Q^{\ell}_{g',[m']}(\delta G)^\nu_{\underline{kr}}	
		\end{equation}
		where $\Q^{\ell}_{g',[m']}(\delta G)_{\underline{kr}}^\nu$ denotes the fixed field of the subgroup $p(\Ker \Phi_{g',[m']}^{\ell}(\delta G)_{\underline{kr}}^\nu)$ as usual. Lemma~\ref{lem:nobrack} then gives $\Q^{\ell}_{g',m'}\subset \Q^{\ell}_{g',m'}(\delta G)_{\underline{kr}}^\nu$.
		
		\subsubsection{} In the rest of this section, we finally establish that $\Q^{\ell}_{g',[m']}(\delta G)_{\underline{kr}}^\nu=\Q^{\ell}_{g,m}(G)_{\underline{kr}}^\nu$ in the case where $X\rightarrow Y$ is a finite étale\footnote{I.e. $\underline{kr}$ is of étale type, that is $\underline{kr}^{et}=\underline{kr}$, see Section~\ref{sub:kret} for definition.} geometric cover and where $G\simeq \Zln$. The finite étale condition guarantees that the inclusion $\iota\colon \pi_1(X)\rightarrow \pi_1(Y)$ induces an inclusion at the pro-$\ell$ completion level $\iota^{\ell}\colon\pi_1^{\ell}(X)\rightarrow \pi_1^{\ell}(Y)$. 
		
		\medskip
		
		Denoting by $\Aut \pi_1^{\ell} (Y)^{X}$ the subgroup of the automorphisms of $\pi_1^{\ell}(Y)$ that stabilizes $\pi_1^{\ell}(X)$, we thus obtain a big commutative diagram
		\[\begin{tikzcd}[row sep = .5cm, column sep=-1em] 
			& \pi_1(X) \arrow [dl] \arrow [rr] \arrow [ddd,hook] & & \pi_1( \mC_{g,[m]}(\Zln)_{\kr}^\nu) \arrow [dl, dashed] \arrow [ddd] \arrow [rr] & & \pi_1(\Mgm(\Zln)_{\kr}^{\nu}) \arrow[dl] \arrow[ddd]\\
			\operatorname{Inn} \pi_1^{\ell}(X) \arrow [rr, crossing over] \arrow[ddd, hook] & & \Aut \pi_1^{\ell}(X) \arrow [rr, crossing over] & & \Out \pi_1^{\ell}(X) \\
			& & \Aut \pi_1^{\ell} (Y)^{X} \arrow[u]\arrow[rr, crossing over]& & \Aut \pi_1^{\ell} (Y)^{X}/\iota^{\ell}(\Inn \pi_1^{\ell}(X)) \arrow[u]  \\
			& \pi_1(Y)\arrow [dl] \arrow [rr] & & \pi_1(\mM_{g',[m']+1}) \arrow [dl] \arrow[rr] & & \pi_1(\mM_{g',[m']}) \arrow[dl]\\
			\operatorname{Inn} \pi_1^{\ell}(Y) \arrow [rr] & & \Aut \pi_1^{\ell}(Y) \arrow[rr] \arrow[from=uu,hook, crossing over]& & \Out \pi_1^{\ell}(Y) \arrow[from=uu, crossing over]\\
		\end{tikzcd}\]
		By tracking the conjugation action of $\pi_1(\mC_{g,[m]}(\Zln)_{\kr}^\nu)$ on $\pi_1^{\ell}(X)$ on the first square of the back face, we see that the dashed arrow $\pi_1(\mC_{g,[m]}(\Zln)_{\kr}^\nu) \rightarrow \Aut \pi_1^{\ell}(X)$ factors by $\Aut \pi_1^{\ell} (Y)^{X}$ through its conjugation action on $\pi_1^{\ell}(Y)$ and the restriction map. 
		
		\begin{theo}\label{theo:fieldQuot} For $(g,m)$ of hyperbolic type, and $\kr$ an abstract Hurwitz data of étale type associated to $\Zln$ with quotient topological data $(g',m')$, we have the following inclusions of $\ell$-monodromy fixed fields
			\[
			\Q^{\ell}_{g',[m']} \subset \Q^{\ell}_{g,[m]}(\Zln)_{\kr}^\nu.
			\]
		\end{theo}
		
		\begin{proof}
			Since $\Q^{\ell}_{g',[m']}\subset \Q^{\ell}_{g',[m']}(\delta G)_{\underline{kr}}^\nu$ by Eq.~\eqref{Eq:gmDetlaFields} it suffices to show the equality $\Q^{\ell}_{g',[m']}(\delta \Zln)_{\underline{kr}}^\nu=\Q^{\ell}_{g,[m]}(\Zln)_{\underline{kr}}^\nu$. We do so by introducing some intermediate fields as can be seen in Diag.~\eqref{Diag:fieldInclBig}.
			
			\medskip
			
			We first have a map
			\[
			\Psi\colon \pi_1^{\ell}(\mC_{g,[m]}(\Zln)_{\kr}^\nu)\longrightarrow \Aut \pi_1^{\ell} (X) \times \Aut \pi_1^{\ell} (Y)^{X}
			\]
			such that $\Phi_{g',[m']}^{\ell}(\delta \Zln)_{\underline{kr}}^{\nu}$ and $\Phi_{g,[m]}^{\ell}(\Zln)_{\underline{kr}}^{\nu}$ are obtained by composing $\Psi$ with the projections and quotients by the inner automorphisms. One checks directly that $\Inn \iota^{\ell} \pi_1^{\ell}(X)$ is a normal subgroup of $\Aut \pi_1^{\ell} (Y)^{X}$. We thus have a quotient map
			\[
			\Aut \pi_1^{\ell} (X) \times \Aut \pi_1^{\ell} (Y)^{X} \longrightarrow \Out \pi_1^{\ell} (X) \times \Aut \pi_1^{\ell} (Y)^{X}/\Inn \iota^{\ell}( \pi_1^{\ell} (X))
			\]
			which by composition with $\Psi$ results in a map
			\[
			S^{\ell}\colon\pi_1(\mM_{g,[m]}(\Zln)^{\nu}_{\underline{kr}})\longrightarrow   \Out \pi_1^{\ell} (X) \times \Aut \pi_1^{\ell} (Y)^{X}/\Inn \iota^{\ell}( \pi_1^{\ell} (X)).
			\]
			
			Considering the quotient map $p_{Y}\colon\Aut \pi_1^{\ell} (Y)^{X}/\Inn \iota^{\ell}( \pi_1^{\ell} (X)) \rightarrow \Out \pi_1^{\ell} (Y)$ and the canonical projections $p_i$, $i=1,2$, of the product $\Out \pi_1^{\ell} (X) \times \Aut \pi_1^{\ell} (Y)^{X}/\Inn \iota^{\ell}( \pi_1^{\ell} (X))$, one observes that by construction
			\[
			\Phi_{g,[m]}^{\ell}(\Zln)_{\underline{kr}}^\nu=p_1 \circ S^{\ell} \text{ and } \Phi_{g',[m']}^{\ell}(\delta \Zln)_{\underline{kr}}^\nu=p_{Y} \circ p_2 \circ S^{\ell}.
			\]
			
			By setting $\Q_S^{\ell}$ to be the fixed field of $p(\Ker S^{\ell})$ and $\Q_{S_Y}^{\ell}$ to be the fixed field of $p(\Ker p_2\circ S^{\ell})$, we obtain the following diagram of inclusions of $\ell$-monodromy fixed fields
			\begin{equation}\label{Diag:fieldInclBig}
				\begin{tikzcd}[row sep=4pt, column sep=6pt]
					\Q_{g,m}^{\ell} \arrow[rr, hook] &   & \Q_{g,[m]}^{\ell}(\Zln)_{\underline{kr}}^\nu  \arrow[drr, hook] &   & \\
					& & & & \Q_{S}^{\ell}\\[6pt]
					\Q_{g',[m']}^{\ell} \arrow[rr, hook] &  &  \Q_{g',[m']}^{\ell}(\delta \Zln)_{\underline{kr}}^\nu \arrow[r, hook] &  \Q_{S_Y}^{\ell} \arrow[ru, hook] &
				\end{tikzcd}		    
			\end{equation}

			It remains to show some equalities. First, $\Q_{S_Y}^{\ell}=\Q_{S}^{\ell}=\Q_{g,[m]}^{\ell}(\Zln)_{\underline{kr}}^\nu$ since by the inclusion $\iota^{\ell}$ we have that $p_2$ restricted to the image of $S^{\ell}$ is injective, and by slimness of $\pi_1^{\ell}(X)$, see Section~\ref{subsub:slimness} for a definition, we have that $p_1$ restricted to the image of $S^{\ell}$ is also injective.
			
			\medskip
			
			For the remaining equality $\Q_{S_Y}^{\ell}=\Q_{g',[m']}^{\ell}(\delta \Zln)_{\underline{kr}}^\nu$, we consider the stack inertia injection $G\subset \mathcal{I}_{\mM,\overline{x}}\hookrightarrow \pi_1(\MgmGkr^{\nu})$ as in \cite{No04}, where $\overline{x}\in \MgmGkr^{\nu}( \overline{K})$ corresponds to the curve $\overline{X}$, and where the injectivity follows from ibid. Theorem~6.2 with the arguments of Remark~4.4 of \cite{CM18}. The injection $G\hookrightarrow \pi_1(\MgmGkr^{\nu})$ can be shown to be independent of the choice of point $\overline{x}$ and maps, through our construction, $G=\Zln$ isomorphically to the quotient $\Inn \pi_1^{\ell}(Y)/\Inn \iota^{\ell}( \pi_1^{\ell} (X)$. Let $\sigma\in p(\Ker p_{Y}\circ p_2\circ S^{\ell})$, which lifts to $\tau \in \pi_1(\MgmGkr^{\nu})$ by definition which in turn maps to $h\in\Inn \pi_1^{\ell}(Y)/\Inn \iota^{\ell}( \pi_1^{\ell} (X))\simeq G$. The element $h^{-1}\tau\in \pi_1(\MgmGkr^{\nu})$ is in $\Ker p_2\circ S^{\ell}$ and verifies $p(h^{-1}\tau)=p(\tau)$. Thus, we have proven that $p(\Ker \Phi^{\ell}_{g',[m']}(\delta \Zln)^\nu) \subset p(\Ker S_Y^{\ell})$, and the reverse inclusion is given by Diag. \ref{Diag:fieldInclBig}.  
		\end{proof}
		
		\medskip
		
		By Theorem~A of \cite{Nacoupl}, Theorem~4.3 of \cite{matsu} and Theorem~3.6 \cite{Tak12} there is an inclusion $\Q^{\ell}_{0,3}\subset \Q^{\ell}_{g,m}$ for all $(g,m)$ of hyperbolic type. Thus we can complete the diagram Diag.~\eqref{Diag:fieldInclBig} of field inclusions as follows.
		
		\begin{cor} \label{cor:inclusiondiagramfields}
			For $(g,m)$ of hyperbolic type and $\kr$ compatible Hurwitz data, we have a diagram of inclusions of fields
			\[\begin{tikzcd}
				\Q^{\ell}_{g,m} \arrow[r, hook] &  \Q^{\ell}_{g,[m]}(\Zln)_{\underline{kr}} \arrow[r, hook] &[2pt] \Q^{\ell}_{g,[m]}(\Zln)_{\underline{kr}}^\nu\\
				\Q^{\ell}_{0,3} \arrow[u, hook,dashed] \arrow[r, hook,dashed] & \Q^{\ell}_{g',m'} \arrow[ur, hook] &  
				&
			\end{tikzcd}\]
		\end{cor}
		
		\begin{rem}\label{rem:beyondOdanoWeight}\mbox{}
			\begin{enumerate}
				\item While for some well-chosen Hurwitz data $\kr$ we have $(g',m')=(0,3)$ in the diagram above, the above references \cite{Nacoupl, matsu, Tak12} are still required for the final comparison of monodromy fields. 
				\item \label{it:rem:OdaWeightProof} In \emph{the setting of Oda's weight conjecture}, where the pro-$\ell$-fundamental groups are replaced by quotients $\pi^{\ell}_1(-)[w]$ with respect to a certain weight filtration $\pi^{\ell}_1(-)(w)$, the map $\pi^{\ell}_1(X)[w]\rightarrow \pi_1^{\ell}(Y)[w]$ fails to be injective. Thus the end of the proof of Theorem~\ref{theo:fieldQuot} does not adapt well, since we can not recover the equality $\Q_{S_Y}^{\ell}(w)=\Q_{S}^{\ell}(w)$, where $\Q_{S_Y}^{\ell}(w)$ and $\Q_{S}^{\ell}(w)$ are defined in the obvious manner. See also Remark~\ref{rem:OdaWeightOK}
			\end{enumerate}
		\end{rem}
		
		\medskip
		
		Establishing the $G$-special loci Oda's conjecture in the case of $G=\Zl$ -- that is that $\Q^{\ell}_{g,m}(\Zl)_{\underline{kr}}$ is independent of the topological and Hurwitz data and indeed equal to $\Q^{\ell}_{0,3}$ -- is thus reduced to establishing the last inclusion $\Q^{\ell}_{g,[m]}(\Zl)_{\underline{kr}}^\nu \subset \Q^{\ell}_{0,3}$. We proceed to do so in the rest of this paper by developing for $G$-special loci a refinement of Ihara-Nakamura's degeneration method used in their original proof of the containment $\Q^{\ell}_{g,m}\subset \Q^{\ell}_{0,3}$ in \cite{iharanak}.

		\section{Maximal degeneration families for $G$-stable compactification} \label{sec:construction}
		After some brief reminder on Deligne's tangential base point on $\mM_{0,4}$, we construct, following \cite{iharanak} for generic curves, some tangential base points on $\Mgm(G)$ as $1$-parameter deformation families $X/\Spf K[[q]]$ of some maximally degenerated $G$-stable curves in some well-chosen strata of $\compMGkr$. These curves are defined as certain $\Zl$-stable $C_r$-diagrams $X^0$ that are obtained, via Grothendieck's formal patching technique, from well-chosen arrangements of so-called Matsumoto-Seyama curves $C_r$. In particular, the associated $\Zl$-quotient curves and their deformation will be the $\PP^1\setminus\{0,1,\infty\}$-diagrams and their canonical $1$-dimensional deformation constructed by Ihara and Nakamura in \cite{iharanak}~2.1.3.
		
		\medskip
		
		We enunciate, under the anabelian slimness hypothesis, some immediate results for the kernel of universal monodromy representations, and for $\Q^{\ell}_{C_r'}$. Consequences for the $\ell$-monodromy fixed fields $\Q^{\ell}_{g,m}$, $\Q^{\ell}_{g,[m]}(\Zl)_{\kr}$ and $\Q^{\ell}_{0,3}$, and for Oda's conjecture are exploited in Section~\ref{sec:conclusion}.
		
		\subsection{Tangential Galois actions and universal monodromy properties} 
		
		\subsubsection{}
		We follow the elementary definition of tangential base point of the survey \cite{Natang} Section~I, that is, for X  connected smooth curve over a field $K$ \emph{a $K$-tangential base point $v$ on $X$} is a morphism $v\colon\Spec K((t))\rightarrow X$ (see ibid. Definition~1.1). 
		
		\medskip
		
		The key feature of such a choice of a $K$-rational tangential base point is, via the field of Puiseux series $\overline{K}\{\{t\}\}$, to provide at once a geometric base point for the étale fundamental group of $X$ and a section to the related homotopy exact sequence:
		
		\begin{equation}\label{eq:FAGexseq}
			\begin{tikzcd}
				1 \arrow[r] & \pi_1(X_{\overline{K}},\vec{v}) \arrow[r]  & \pi_1(X,\vec{v}) \arrow[r] & G_K \arrow[r] \arrow[l, "s_{v}", bend right = 50] & 1 \\
			\end{tikzcd}
		\end{equation}
		
		In other words, one obtains a specific $G_K$-action $\varphi_{\vec{v}}$ on $\pi_1(X_{\overline{K}},v)$ given by conjugation which lifts the canonical outer Galois action
		\begin{equation}
			\begin{tikzcd}
				G_K\arrow[r,"\varphi_{v}"] \arrow[rd, "\varphi_X"']& \Aut \pi_1(X_{\overline{K}}, v) \arrow[d]\\
				& \Out \pi_1(X_{\overline{K}},v)
			\end{tikzcd}
		\end{equation}
		and can be chosen to reflect some good arithmetic properties of $X$. More explicitly, the $G_K$ action $\varphi_v$ is given, via the function fields of $X$, by the action on the coefficients of the formal series in $\overline{K}\{\{t\}\}$, see also Eq.~\eqref{it:GMtgtAction} below.
		
		
		\pagebreak[4]
		
		\begin{rem}\mbox{}
			\begin{enumerate}			
				\item By the valuative criterion of properness this is equivalent to giving a map $v\colon \Spec K[[t]]\rightarrow \overline{X}$ where $\overline{X}$ is the compactification of $X$ (i.e., $X\subset \overline{X}$ is a Zariski open given by finitely many punctures of the proper curve $\overline{X}$).
				
				\item \label{it2:rem:GMtgt} By Grothendieck-Murre theory, the category $\operatorname{Rev}^D(\overline{X})$ of finite étale coverings of $\overline{X}$ tamely ramified along the divisor $D=\overline{X}\setminus X$ is equivalent to the category of finite étale coverings of $X$. The choice of a tangential base point gives a fiber functor of this Galois category in the following way. Let $Y\in \operatorname{Rev}^D(\overline{X})$ and $B$ the $K[[t]]$-algebra obtained by the pullback of $Y$ along our tangential base point $v$. With this formalism, the fiber functor $\vec{v}$ is defined by
				\begin{equation}\label{it:GMtgtAction} 
					\begin{array}{cccc}
						\vec{v}\colon &  \operatorname{Rev}^D(\overline{X})   & \longrightarrow  & \mathrm{Set}  \\
						& Y    & \longmapsto & \Hom_{K[[T]]}(B, \overline{K}\{\{t\}\}). 
					\end{array}
				\end{equation}
				
				\item 
				The above formalism provides a fundamental group $\pi_1^{D}(\overline{X},v)$ which is canonically isomorphic to $\pi_1(X,v)$ and carries the same tangential Galois action.
			\end{enumerate}
		\end{rem}
		
		\medskip
		
		For $X=\PP_{\Q}^1\setminus \{ 0,1, \infty\}$ let us denote the set of fiber functor associated, as in \ref{it2:rem:GMtgt} above, to Deligne-Ihara's original $\Q$-tangential base points by
		\[
		\mathbb{B}=\{\vec{01}, \vec{0\infty}, \vec{10}, \vec{1\infty}, \vec{\infty 1}, \vec{\infty 0}\}
		\]
		where for example $\vec{01}\colon \Spec \Q((t))\to \PP_{\Q}^1\setminus \{ 0,1, \infty\}$ and $\vec{0\infty}\colon \Spec \Q((-t))\to \PP_{\Q}^1\setminus \{ 0,1, \infty\}$, and refer to the Appendix of \cite{IHA93} for further details on the associated $G_\Q$-action. For our study, the main property of these tangential base points is that
		\begin{equation}\label{eqDelTgt}
			\Ker \varphi^{\ell}_{\vec{ij}}= \Ker \varphi^{\ell}_{\PP^1\setminus \{ 0,1, \infty\}}\quad \text{for every } \vec{ij}\in \mathbb{B}.
		\end{equation}
		
		\medskip
		
		While even the simplest rational scaling of the parameter, see for example $\vec{01}$ vs $\vec{0\infty}$ above or \cite{tsuno} Section~1.5, changes the tangential Galois action, we have the following Galois invariance property.

		\begin{lem} \label{lem:basepointaction}
			The $G_K$-action induced by a $K$-rational tangential base point $v\colon \Spec K((t))\rightarrow X$ depends only on the closed point $x\in \overline{X}(K)$ in the closure of the image of $v$ and the class of the image of $t$ in the cotangent space $\mfm_x/\mfm_x^2$.
		\end{lem}
		
		\begin{proof}
			Let $x\in \overline{X}$ be a closed $K$-rational point. It suffices to show that if $t$ and $t'$ are both uniformizers at $x$ (i.e., we have $\widehat{\OO}_{\overline{X},x}\simeq K[[t]] \simeq K[[t']]$ and $t'=t(1+tF)$ in $K[[t]]$ with $F\in K[[t]]$) then the isomorphism $\delta_{t',t}\colon \overline{K}\{\{t'\}\} \rightarrow \overline{K}\{\{t\}\}$ is $G_K$-equivariant. But as $\delta_{t',t}$ is defined by ${t'}^{\frac{1}{N}} \mapsto t^{\frac{1}{N}}(1+tF)^{\frac{1}{N}}$ for $N\geq 0$ this comes from the fact that $(1+tF)^{\frac{1}{N}}=G_N$ with $G_N\in K[[t]]$ by the series expansion of $(1+tF)^{\frac{1}{N}}$. 
			
			\medskip
			
			Indeed, let $v_t$ (resp.~$v_t'$) be the tangential base points given by $t$ (resp.~$t'$) and denote by $\varphi_{v_t}$ (resp.~$\varphi_{v_t'}$) the associated tangential $G_K$-action. Let $\sigma\in G_K$ and $f=\sum\limits_{k} a_k {t'}^{\frac{k}{N}}\in \mM_{v_{t'}}\subset K\{\{t'\}\}$. Then we have
			\begin{align*}
				\sigma^{-1}_{v_t}\circ \delta_{t',t} \circ \sigma_{v_{t'}}(f) &= \sigma^{-1}_{v_t} \circ \delta_{t',t}( \sum\limits_{k} \sigma(a_k) {t'}^{\frac{k}{N}}) \\
				&= \sigma^{-1}_{v_t}( \sum\limits_{k} \sigma(a_k) {t}^{\frac{k}{N}}G_N) \\
				&=\sum\limits_{k}a_k {t}^{\frac{k}{N}}G_N 
			\end{align*}
			that is
			\[
			\sigma^{-1}_{v_t}\circ \delta_{t',t} \circ \sigma_{v_{t'}}(f) = \delta_{t',t}(f)
			\]
			which shows that $\sigma_{v_t}^{-1}\circ \delta_{t',t}\circ \sigma_{v_{t'}}= \delta_{t',t}$ and thus $\varphi_{v_t'}=\varphi_{v_t}\circ \delta_{t',t}$ as intended.
		\end{proof}
		
		\subsubsection{} \label{subsub:slimness}
		
		We recall that \emph{a profinite group is said to be slim} if any of its open subgroup has trivial  centralizer. Examples of slim groups include the absolute Galois group of rational numbers and the pro-$\ell$ fundamental group of hyperbolic curves, see \cite{MOCTAM08} Proposition~1.4.

		\medskip
		
		We record the following inclusions between the $\ell$-monodromy fixed fields of the various tangential and non-tangential Galois actions in the case of étale coverings.
		
		\begin{lem}\label{lem:inclusion-étale}
			Let $\psi\colon X\rightarrow Y$ be a finite étale covering of geometrically irreducible curves over a field $K$ of degree a power of $\ell$. Let $v\colon \Spec K((t)) \rightarrow X$ be a tangential base point on $X$ and $\psi(v)$ the induced tangential base point on $Y$. We have the following inclusions of subgroups of $G_K$:
			\begin{enumerate}
				\item\label{it1:lem:inclusion-étale} $\Ker \varphi^{\ell}_{\vec{\psi(v)}} \subset \Ker \varphi^{\ell}_{\vec{v}}$
				\item\label{it2:lem:inclusion-étale} $\Ker \varphi^{\ell}_{\vec{v}} \subset \Ker \varphi^{\ell}_X$ and $\Ker \varphi^{\ell}_{\vec{\psi(v)}} \subset \Ker \varphi^{\ell}_Y.$
			\end{enumerate}
			Furthermore, when $\pi_1^\ell(Y,\vec{\psi(v)})$ is slim we have $\Ker \varphi_{\vec{v}}= \Ker \varphi_{\vec{\psi(v)}}$ and $\Ker \varphi^{\ell}_X \subset \Ker \varphi^{\ell}_Y$.
		\end{lem}
		
		\begin{proof}
			The homotopy exact sequence for $X$ and $Y$ and the covering map $\psi$ gives the diagram
			\[\begin{tikzcd}
				1 \arrow[r] & \pi_1(X_{\overline{K}},\vec{v}) \arrow[r] \arrow[d, hook]  & \pi_1(X,\vec{v}) \arrow[r] \arrow[d] & G_K \arrow[r] \arrow[d, equal] \arrow[l, "s_{v}", bend right = 50] & 1 \\
				1 \arrow[r] & \pi_1(Y_{\overline{K}},\vec{\psi(v)}) \arrow[r]  & \pi_1(Y,\vec{\psi(v)}) \arrow[r]  & G_K \arrow[r] \arrow[l, "s_{\psi(v)}"', bend right = 50] & 1 \\
			\end{tikzcd}
			\]
			that is commutative by definition of $\psi(v)$ and the étaleness of $\psi$.
			One thus recovers, via $\pi_1(Y,\vec{\psi(v)}) \to \Aut \pi_1^{\ell}(Y_{\overline{K}},\vec{\psi(v)})$ whose image stabilizes $\pi_1^\ell(X_{\overline{K}},\vec{v})$, the monodromy action $\varphi^{\ell}_{\vec{v}}$ as the composition 
			\[
			G_K\rightarrow \Aut \pi_1^{\ell}(Y_{\overline{K}},\vec{\psi(v)})^{X_{\overline{K}}} \rightarrow  \Aut \pi_1^{\ell}(X_{\overline{K}},\vec{v})
			\]
			which leads to the inclusion given in \ref{it1:lem:inclusion-étale}. In the case of slimness the right restriction map is injective, which yields the equality. 
			
			\medskip
			
			The remaining inclusions are obtained by adding the following commutative diagram
			\[\begin{tikzcd}
				\Aut \pi_1^{\ell}(X_{\overline{K}},\vec{v}) \arrow[d] & \Aut \pi_1^{\ell}(Y_{\overline{K}},\vec{\psi(v)})^{X_{\overline{K}}} \arrow[l] \arrow[d, "d_X"]  \\
				\Out \pi_1^{\ell}(X_{\overline{K}},\vec{v}) & \Aut \pi_1^{\ell}(Y_{\overline{K}},\vec{\psi(v)})^{X_{\overline{K}}}/ \Inn \pi_1^{\ell}(X_{\overline{K}},\vec{v}) \arrow[l] \arrow[d, "d_Y"] \\
				& \Out \pi_1^{\ell}(Y_{\overline{K}},\vec{\psi(v)})\\
			\end{tikzcd}
			\]
			
			The inclusions of \ref{it2:lem:inclusion-étale} are thus direct by the diagram and the definitions of the maps involved. To see the remaining inclusion, we remark that by slimness $\Ker \varphi^{\ell}_X= \Ker d_X\circ \varphi^{\ell}_{\vec{\psi(v)}}$, and the inclusion follows as $\varphi^{\ell}_Y= d_Y\circ d_X\circ \varphi^{\ell}_{\vec{\psi(v)}}$.
		\end{proof}
		
		\subsection{The Matsumoto-Seyama curves}\label{sec:seyamacurves}
		We now introduce the Matsumoto-Seyama curves $C_r$, for $r\in \{0,1,\dots, \ell-2\}$, that live in certain special loci $\Mgm(\Zl)_{\kr}$ and that have $\PP^1_K$ as $\Zl$-quotient,  where $K$ denotes $\Q(\mu_{\ell})$. 
		
		\medskip
		
		For $r\in \{1,\dots, \ell-2\}$, the curves $C_r$ are those of \cite{seya}, that is, some smooth projective curves of genus $g={\ell-1}/{2}$ that are birationally equivalent to the affine curve
		\begin{equation}
			y^r(y-1)=x^{\ell}\text{ with } \Zl\text{-action }\begin{cases}
				\text{given by } x\mapsto \zeta_{\ell} x\\
				\text{ramified at } P_{r,0}, P_{r,1}, P_{r,\infty}\text{ over } 0, 1,\infty.
			\end{cases}
		\end{equation} 
		The quotient $\psi\colon C_r\setminus \{P_{r,0}, P_{r,1}, P_{r,\infty}\} \rightarrow \PP^1_K\setminus \{0,1,\infty\}$ is finite étale and Galois of group $\Zl$. The abstract Hurwitz data of $C_r$ is $\underline{k}=(r,1,-(r+1))$ which, when $r$ varies, is seen to represent every possible abstract Hurwitz data of a $\Zl$-curve with three ramified points. 
		
		\subsubsection{}
		The set of curves $\{C_r \mid r=1,\dots, \ell-2\}$ admits an $\mfS_3$-action that is compatible with the $\Zl$-action and, in particular, with the $\mfS_3$-action on $\PP^1_K\setminus\{0,1,\infty\}$ through the quotient map, see \cite{seya} Corollary~2.5. This allows us to define, for every $r$, the tangential base points on $C_r'=C_r\setminus \{P_{r,0}, P_{r,1}, P_{r,\infty}\} $ at the punctures by doing so at $P_{r,1}$ . 
		
		\medskip
		
		Indeed, for $\sigma\in \mfS_3$ we have $\sigma(P_{r,1})=P_{\sigma(r),\sigma(1)}$ so that, for every $r\in \{1,\dots, \ell-1\}$ and every $P\in \{P_{r,0}, P_{r,1}, P_{r,\infty}\}$, there is an element $\sigma\in \mfS_3$ such that $P$ is the image of $P_{r,1}$ for some $r$. Now, the smooth affine open $U=C_r\setminus \{P_{r,0}, P_{r,\infty}\}$ is given by  $U=\Spec K[x,y, \frac{1}{y}]$ where $x^{\ell}=y^r(y-1)$. Looking at the equation we see that $x$ is a uniformizer at $P_{r,1}$ and we have $\psi(x)=x^{\ell}=y^r(y-1)$ where $\psi$ is the quotient map to $\PP^1_K$. 
		
		\begin{lem}\label{lem:TcompDelAct}
			The tangential base point $T^r_{10}\colon \Spec K((t)) \rightarrow  C_r'$ defined by $t \mapsto \zeta_{2\ell} x $ induces a tangential base point $\psi(T^r_{10})$ on $\PP^1_K\setminus\{0,1,\infty\}$ that defines the same $G_K$-action on $\pi_1(\PP^1_{\overline{K}}\setminus\{0,1,\infty\},\vec{10})$ as $\overrightarrow{10}$.
		\end{lem}
		
		\begin{proof}
			By Lemma \ref{lem:basepointaction} it suffices to check that $\psi(T^r_{10})$ and $\overrightarrow{10}$ have, after taking the closure, the same closed points in $\PP^1_K$ and the same class in $\mfm_1/\mfm_1^2.$ The first part is obvious. For the second one, by definition, we have that $\psi(T^r_{10})$ is $-y^r(y-1)\in K[[y-1]]\simeq \widehat{\OO}_{\PP^1_K,1}$ so that its class modulo $\mfm_1^2$ is equal to $-1$ as required.  
		\end{proof}
		
		As stated before, by using the $\mfS_3$-action on the previous subset of Matsumoto-Seyama curves, we obtain tangential base points $T^r_{ij}$ for $i,j\in \{0,1,\infty\}$, whose set of associated fiber functors on the categories of finite étale covers $\mathrm{Et}(C_r')$ we denote by
		\[
		\B^r=\{\overrightarrow{T^r_{ij}} \mid i,j\in \{0,1,\infty\}\},\text{ for } r\in \{1,\dots, \ell-2\}.
		\]
		These tangential base points induces same tangential $G_K$-actions on the fundamental group of $\PP^1_{\overline{K}} \setminus\{0,1,\infty\}$ given by Deligne-Ihara as in Lemma~\ref{lem:TcompDelAct}. 
		
		\begin{theo} \label{the:actioncruvep1}
			The $G_K$-action defined by the $\overrightarrow{T^r_{ij}}$s on the groupoid $\Pi_1({C_r'}_{,\overline{K}}, \mathbb{B}^r)$ induces a $G_K$-action on the groupoid $\Pi_1(\PP_{\overline{K}}^1\setminus\{0,1,\infty\}, \mathbb{B})$ that is compatible with the Deligne-Ihara one. Furthermore, an element of $G_K$ acts trivially on $\Pi_1^{\ell}({C_{r}'}_{,\overline{K}}, \mathbb{B}^r)$ if and only if it acts trivially on $\Pi_1^{\ell}(\PP_{\overline{K}}^1\setminus\{0,1,\infty\}, \mathbb{B})$.
		\end{theo}
		
		\begin{proof}
			The first part of the statement is the result of the previous paragraph. For the second part, let $\sigma\in G_K$. As the tangential base points of $\mathbb{B}^r$ are $K$-rational, the action of $\sigma$ on $\Pi_1({C_r'}_{\overline{K}}, \mathbb{B}^r)$ stabilizes each fundamental group or set of étale paths. Now as this action is compatible with the one on $\Pi_1(\PP_{\overline{K}}^1\setminus\{0,1,\infty\}, \mathbb{B})$ and each of the inclusions maps between $\Pi_1({C_r'}_{,\overline{K}}, \overrightarrow{T_{ij}}, \overrightarrow{T_{jk}})$ and $\Pi_1(\PP_{\overline{K}}^1\setminus\{0,1,\infty\}, \vec{ij}, \vec{jk})$ remains injective after passing to the pro-$\ell$-completion for all $i,j,k\in \{0,1,\infty\}$, it follows that the reverse implication holds. By Lemma~\ref{lem:inclusion-étale}, it also holds that $\Ker \varphi^{\ell}_{\overrightarrow{T_{ij}}} = \Ker \varphi^{\ell}_{\vec{ij}}$ for all $i,j\in \{0,1,\infty\}$. Thus, if $\sigma$ acts trivially on  $\Pi_1(\PP_{\overline{K}}^1\setminus\{0,1,\infty\}, \mathbb{B})$, it acts trivially on each of the fundamental groups appearing in $\Pi_1({C_r'}_{\overline{K}}, \mathbb{B}^r)$, and thus on the whole groupoid. 
		\end{proof}
		
		In what follows $r$ will be omitted from notations when clear from context.
		
		\subsubsection{}
		For $r=0$, we consider the covering of $\PP^1_K$ given by
		\[
		C_0\colon x=y^{\ell} \text{ with usual }\Zl\text{-action having two ramified points } 0 \text{ and }\infty
		\]
		with abstract Hurwitz data $\underline{k}=(1,-1)$. The $\ell+2$-marking is given by the two ramified points and by the unramified points $P_1,\dots, P_{\ell}$ of the fiber at $1$. We further set 
		\[
		\B^0=\{\overrightarrow{T^0_{0\infty}}, \overrightarrow{T^0_{\infty0}}\} \text{ and } C_0'=C_0\setminus \{0,\infty, P_1,\dots, P_{\ell}\},
		\]
		where the fiber functors $\overrightarrow{T^0_{0\infty}}$ and $\overrightarrow{T^0_{\infty0}}$ are induced by the tangential base points associated to the parameter $x$ and $\frac{-1}{x}$, and which are direct lifts of the Deligne tangential base points $\overrightarrow{0\infty}$ and $\overrightarrow{\infty0}$.
		
		\begin{prop}
			The action of $G_K$ on $\Pi_1({C_0'}_{\overline{K}}, \B^0)$ is compatible with its action on $\Pi_1(\PP_{\overline{K}}^1\setminus\{0,1,\infty\}, \mathbb{B})$. Furthermore, an element of $G_K$ acts trivially on $\Pi_1^{\ell}({C_0'}_{\overline{K}}, \B^0)$ if and only if it acts trivially on $\Pi_1^{\ell}(\PP_{\overline{K}}^1\setminus\{0,1,\infty\}, \mathbb{B}).$
		\end{prop}
		
		\begin{proof}
			The only part of the statement that is not already proven is a direct consequence of \cite{iharanak} Corollary~4.1.4~(ii). 
		\end{proof}
		
		\subsubsection{}
		We finish this section by showing that the $\ell$-monodromy fixed field of the Matsumoto-Seyama curves is $\Q^{\ell}_{0,3}$.
		
		\begin{cor}\label{cor:seyamagenus0field}
			We have $\Q^{\ell}_{C_r'}= \Q^{\ell}_{0,3}$ for all $r\in \{0,\dots, \ell-2\}$.
		\end{cor}
		
		\begin{proof}
			For $r\in\{0,\dots, \ell-2\}$ given, it follows from Lemma~\ref{lem:inclusion-étale} that $\Ker \varphi^{\ell}_{\overrightarrow{T_{0\infty}}}= \Ker \varphi^{\ell}_{\vec{01}}$ as $\pi_1(\PP^1_K\setminus \{0,1,\infty\}, \vec{0\infty})$ is slim. From the same lemma, we also get the inclusions
			\[
			\Ker \varphi^{\ell}_{\overrightarrow{T_{0\infty}}} \subset \Ker \varphi^{\ell}_{C_r'} \subset \Ker \varphi^{\ell}_{\PP^1\setminus\{0,1,\infty\}}.
			\]
			Since the two outmost terms are equal as in Eq.~\eqref{eqDelTgt}, it follows that $\Ker \varphi^{\ell}_{C_r'}=\Ker \varphi^{\ell}_{\PP^1\setminus\{0,1,\infty\}}$, thus the desired equality.
		\end{proof}

		\begin{rem}\mbox{}
			\begin{enumerate}
				\item At this stage, one can already obtain, by following Matsumoto's approach as in \cite{matsu}, that $\Q^{\ell}_{g,m}=\Q^{\ell}_{0,3}$ for the specific values of $(g,m)=((\ell-1)/2),3)$ and $(g,m)=(0,\ell+2)$.
				
				\item The curves introduced in this section are chosen so that the corresponding stacks $\mM_{g,[m]}(\Zl)_{\underline{kr}}$ have only one geometric point.
			\end{enumerate}
		\end{rem}
		
		\subsection{Diagrams in the $\Zl$-stable compactification}
		Similarly to the $\PP^1_K\setminus\{0,1,\infty\}$-diagrams construction of \cite{iharanak}~1.2, we construct some $\Zl$-stable $C_r$-diagram $X^0$ over a field $K$, here as gluing the previously defined Matsumoto-Seyama $\Zl$-curves. 
		
		\subsubsection{}
		While the gluing, or clutching, of marked points for stable curves can be found in details in \cite{Knu83}, the similar gluing for curves with $G$-action requires an additional constraint as follows.
		
		\medskip
		
		Consider two curves $C_r$ and $C_{r'}$ with $r,r'\in \{1,\dots, \ell-2\}.$ The gluing of both curves at the points $P_{r,1}$ and $P_{r',1}$ can be constructed as the union 
		\[
		C_{r,r'}^{1,1}=C_r \times \{P_{r',1}\} \cup C_{r'}\times \{ P_{r,1}\} \text{ in the fiber product } C_r\times_{\Spec K} C_{r'}.
		\]
		The result of the gluing is a curve $X^0$ of genus $\ell-1$ with $2$ irreducible components and $4$ marked points given by $\{P_{r,0}, P_{r,\infty}, P_{r',0}, P_{r',\infty} \}$, that is equipped with a $\Zl$-action by pullback of the action on the product. 
		
		\medskip
		
		For $X^0$ to be a $G$-stable curve, the $G$-actions must be chosen such that Hurwitz data at the points $P_{r,1}$ and $P_{r',1}$ have opposite characters, see \cite{BERO07} Section~4.1, which is easily done by choosing that $G=\Zl$ acts by $x\mapsto \overline{\zeta_{\ell}}x$ on $C_{r'}$ and by $x\mapsto \zeta_{\ell}x$ on $C_r$. The same construction can be made by gluing together any two ramified points $P_{r,i}$ and $P_{r',j}$ into a curve $C_{r,r'}^{i,j}$, where $i,j\in \{0,1,\infty\}$ denotes which points are glued.
		
		\medskip
		
		Note that the $\mfS_3$-action on the curves $(C_r)_{r\in \{1,\dots, \ell-2\}}$ extends naturally to a $\mfS_3\times \mfS_3$-action on the fiber products $(C_r\times_{\Spec K} C_{r'})_{r,r'\in \{1,\dots, \ell-2\}}$ of such curves. One checks that this action stabilizes the closed subsets $(C_{r,r'}^{i,j})_{r,r'\in \{1,\dots\ell-2\},~i,j\in \{0,1,\infty\}}$ globally, that is for $\sigma, \tau\in \mfS_3\times \mfS_3$ we have $(\sigma,\tau)\cdot C_{r,r'}^{i,j}=C_{\sigma(r), \tau(r')}^{\sigma(i), \tau(j)}$. It results that the affine neighborhood of $C_{r,r'}^{i,j}$ with the $4$ marked points removed is, for some $r$, always isomorphic to 
		\[ 
		C_{r,r'}^{1,1}\setminus \{P_{r,0}, P_{r,\infty}, P_{r',0}, P_{r',\infty} \} = \Spec K[x,y,x',y'][\frac{1}{y},\frac{1}{y'}]/(xx') 
		\]
		which serves as a model for the construction of the $U_\mu$s as in Section~\ref{subsub:X0opens}.

		\subsubsection{}
		We will build our $\Zl$-stable $C_r$-diagrams from the two types of Matsumoto-Seyama curves $C_r$ of Section~\ref{sec:seyamacurves}. Recall that the genus $0$ curves have two distinguished rational sections given by the ramified points, and that the genus $(\ell-1)/2$ ones have three.
		
		\begin{defin}\label{def:Gdiagram} A \emph{$\Zl$-stable $C_r$-diagram} is a connected curve $X^0$ over $K$ that is defined by the following data:
			\begin{enumerate}
				\item\label{it:DiagCurves} \emph{A finite collection of curves} $X^0_{\lambda}$ ($\lambda\in \Lambda\sqcup \Lambda'$) where $X^0_{\lambda}$ is either isomorphic to $C_{r}$ with $r\geq 1$ if $\lambda\in \Lambda$ or to $C_{0}$ if $\lambda\in \Lambda'$. 
				
				\item\label{it:DiagSec} \emph{A finite collection of pairs of distinguished section} $P^0_{\mu}$ ($\mu\in M$) of the $X^0_{\lambda}$, $\lambda\in \Lambda\sqcup \Lambda'$. The pairs $P^0_{\mu}$ are such that the Hurwitz data at those sections are opposite and such that two distinct pairs $P^0_{\mu}$ and $P^0_{\mu'}$ ($\mu\neq \mu'$) have no common element. Let $\mu\in M$ and set $\lambda(\mu)=(\lambda,\lambda')$ where the sections of $P^0_{\mu}$ land in $X^0_{\lambda}$ and $X^0_{\lambda'}$. 
			\end{enumerate}
			The curve $X^0$ is obtained from the disjoint union $\bigsqcup\nolimits_{\lambda\in \Lambda} X^0_{\lambda}$ by identifying the pair of points given by the $P^0_{\mu}$. Given a $\Zl$-stable $C_r$-diagram $X^0$ we shall denote by $Q_v^0$, $v\in N$, the distinguished sections of $X^0$ coming from the $X^0_{\lambda}$ that do not appear in the pairs $P^0_{\mu}$, $\mu \in M$.
		\end{defin}
		
		The isomorphisms of \ref{it:DiagCurves} come with choices of variables $x_{\lambda},y_{\lambda}$ and choices of tangential base points $T^{\lambda}_{ij}\colon \Spec K((t)) \rightarrow X^0_{\lambda}$ with the properties of the ones defined in Section \ref{sec:seyamacurves}. The corresponding set of fiber functors will be denoted by $\B^r_{\lambda}$. We will omit $\lambda$ and $r$ from the notations when it is clear from context.
		
		\subsubsection{}
		The following three kinds of $\Zl$-stable $C_r$-diagrams will be used as basic building blocks for the special fiber of our $1$-parameter deformation families. 
		
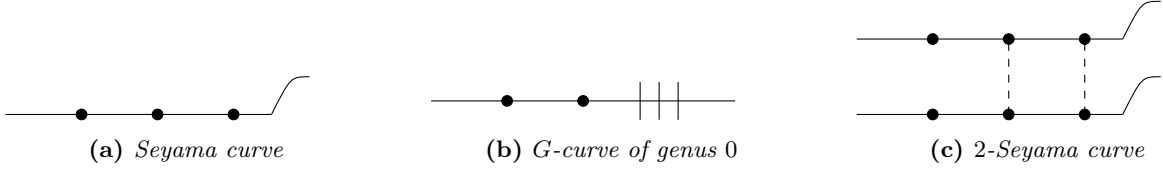
\begin{figure}	
			\begin{subfigure}[b]{0.3\textwidth}
				\begin{tikzpicture} 
					\draw (-2,0) -- (1.5,0);
					\filldraw (0,0) circle (2pt);
					\filldraw (-1,0) circle (2pt);
					\filldraw (1,0) circle (2pt);
					\draw (1.5,0) .. controls (1.75,0.5) .. (2,0.5);
				\end{tikzpicture} 
				\caption{Seyama curve}\label{fig:SeyamCurve}
			\end{subfigure}
			\hfill
			\begin{subfigure}[b]{0.3\textwidth}
				\begin{tikzpicture}
					\draw (-2,0) -- (2,0);
					\filldraw (0,0) circle (2pt);
					\filldraw (-1,0) circle (2pt);
					\draw (0.75,-0.25) -- (0.75,0.25);
					\draw (1,-0.25) -- (1,0.25);
					\draw (1.25,-0.25) -- (1.25,0.25);
				\end{tikzpicture} 
				\caption{$G$-curve of genus $0$}\label{fig:ZeroCurve}
			\end{subfigure}
			\hfill		
			\begin{subfigure}[b]{0.3\textwidth}
				\begin{tikzpicture}
					\draw (-2,0) -- (1.5,0);
					\filldraw (0,0) circle (2pt);
					\filldraw (-1,0) circle (2pt);
					\filldraw (1,0) circle (2pt);
					\draw (1.5,0) .. controls (1.75,0.5) .. (2,0.5);
					
					\draw (-2,-1) -- (1.5,-1);
					\filldraw (0,-1) circle (2pt);
					\filldraw (-1,-1) circle (2pt);
					\filldraw (1,-1) circle (2pt);
					\draw (1.5,-1) .. controls (1.75,-0.5) .. (2,-0.5);
					\draw [dashed] (0,-1) -- (0,0);
					\draw [dashed] (1,-1) -- (1,0);
				\end{tikzpicture} 
				\caption{$2$-Seyama curve}\label{fig:TwoSeyama}
			\end{subfigure}
			\caption{Elementary building blocks of $\Zl$-stable $C_r$-diagrams}\label{fig:BlockCurves}
		\end{figure}

		\begin{itemize}
			\item \textbf{Seyama curve} (Fig.~\ref{fig:SeyamCurve}): a curve of genus $g=(\ell-1)/2$ with $\nu=3$ ramified points and $\underline{k}$ free;
			
			\item \textbf{A $\Zl$-curve of genus $0$} (Fig.~\ref{fig:ZeroCurve}): a curve of genus $g=0$ with $\nu=2$ ramified points, $\ell$ unramified points and $\underline{k}=(1,-1)$;
			
			\item \textbf{A $2$-Seyama curve} (Fig.~\ref{fig:TwoSeyama}): a curve of genus $g=\ell$ with $\nu=2$ ramified points and $\underline{k}=(1,-1)$. These curves are obtained by gluing two Seyama curves twice. The correspond data of the $\Zl$-stable $C_r$-diagram is $\Lambda=\{1,2\}$ where $X^0_{1}$, $X^0_{2}$ are Matsumoto-Seyama curves isomorphic to $C_r$ with $r\in\{1,\dots, \ell-2\}$ and Hurwitz data $\underline{k}=(1,a,-(1+a))$ with ramified points $\{\nu_1,\nu_2,\nu_3\}$ and $\underline{k'}=(-1,-a,1+a)$ with ramified points $\{\nu'_1,\nu'_2,\nu'_3\}$. We have $M=\{1,2\}$ and the pairs $P^0_{1},$ $P^0_{2}$ are respectively $(\nu_1,\nu'_1)$ and $(\nu_2,\nu'_2)$. The remaining set of distinguished sections is $\{\nu_3,\nu'_3\}= \{Q^0_v\}_{v\in N}$.			
		\end{itemize}
		
		\medskip
		
		In Fig.~\ref{fig:BlockCurves} above, the bold points represent ramified points under the $G$-action, the lined markings represent the unramified points ($\ell=3$ here), and the dashed lines represent the glued points. The hook at the end of the genus $g=(\ell-1)/{2}$ curves is to differentiate them from the genus $0$ ones, which are represented by straight lines. 
		
		\medskip
		
		Recall that we denote by $\compMG$ the stable compactification of $\MgmG^\nu$, and accordingly $\CompMgmGkr$ the closure of $\MgmGkr$ in $\compMG$.
		\begin{prop}\label{pro:diagramexistence}
			Let $g,m$ and $\underline{kr}$ be given as below, such that $\MZlkr$ is non-empty, then there exists a $\Zl$-stable $C_r$-diagram in the boundary of $\compMZlkr$.
		\end{prop} 
		
		Let us first recall that, by \cite{CM15} Proposition~3.7, the locus $\mM_{g,[m]}(\Zl)_{\underline{kr}}$ is non-empty as soon as $g$ can be obtained by the Hurwitz formula
		\[
		g=(N-2)\frac{\ell-1}{2}+g'\ell \text{ with } g'\geq 0 \text{ and } N \geq 0,1
		\]
		where $N$ is the number of ramified points in the cover. This is a particular instance of the Frobenius coin problem and it is thus known that all $g\geq (\frac{\ell-1}{2})(\frac{\ell-3}{2})$ are attainable with $N-2\geq 0$, as well as one element of each pair $(k, \ell\frac{\ell-1}{2}-\ell-\frac{\ell-1}{2}-k)$ for $k\in \{ 0,\dots , \frac{\ell-1}{2}\frac{\ell-3}{2}-1\}$. When $g\in \{ 0,\dots , \frac{\ell-1}{2}\frac{\ell-3}{2}-1\}$ is attainable only by the choice $N=0$ we say that \emph{$g$ is an unramified case}. For example, this is the case for $g=1$ by considering the translation action by a choice of order $\ell$ point on an elliptic curve. 
		
\begin{figure}[h]
			\centering
			\begin{subfigure}[c]{0.45\textwidth}
				\centering	
				\resizebox{\columnwidth}{!}{%
					\begin{tikzpicture}
						\draw (6,-2) -- (6,1.5);
						\filldraw (6,0) circle (2pt);
						\filldraw (6,-1) circle (2pt);
						\filldraw (6,1) circle (2pt);
						\draw (6,1.5) .. controls (6.5,1.75) .. (6.5,2);
						
						\draw (5,-2) -- (5,1.5);
						\filldraw (5,0) circle (2pt);
						\filldraw (5,-1) circle (2pt);
						\filldraw (5,1) circle (2pt);
						\draw (5,1.5) .. controls (5.5,1.75) .. (5.5,2);
						\draw [dashed] (5,0) -- (6,0);
						\draw [dashed] (5,1) -- (6,1);

						\draw (4,-2) -- (4,1.5);
						\filldraw (4,0) circle (2pt);
						\filldraw (4,-1) circle (2pt);
						\filldraw (4,1) circle (2pt);
						\draw (4,1.5) .. controls (4.5,1.75) .. (4.5,2);
						
						\draw (3,-2) -- (3,1.5);
						\filldraw (3,0) circle (2pt);
						\filldraw (3,-1) circle (2pt);
						\filldraw (3,1) circle (2pt);
						\draw (3,1.5) .. controls (3.5,1.75) .. (3.5,2);
						\draw [dashed] (3,0) -- (4,0);
						\draw [dashed] (3,1) -- (4,1);
						
						\draw [dashed] (4,-1) -- (5,-1);  
						
						\draw (2,-2) -- (2,1.5);
						\filldraw (2,0) circle (2pt);
						\filldraw (2,-1) circle (2pt);
						\filldraw (2,1) circle (2pt);
						\draw (2,1.5) .. controls (2.5,1.75) .. (2.5,2);
						
						\draw (1,-2) -- (1,1.5);
						\filldraw (1,0) circle (2pt);
						\filldraw (1,-1) circle (2pt);
						\filldraw (1,1) circle (2pt);
						\draw (1,1.5) .. controls (1.5,1.75) .. (1.5,2);
						\draw [dashed] (1,0) -- (2,0);
						\draw [dashed] (1,1) -- (2,1);
						
						\draw [dashed] (2,-1) -- (2.25,-1);
						\draw [dashed] (3,-1) -- (2.75,-1);
						
						\draw [loosely dotted] (2.75,0.5) --  (2.25,0.5);
						
						\draw (0,-2) -- (0,1.5);
						\filldraw (0,0) circle (2pt);
						\filldraw (0,-1) circle (2pt);
						\filldraw (0,1) circle (2pt);
						\draw (0,1.5) .. controls (0.5,1.75) .. (0.5,2);
						
						\draw [dashed] (0,1) -- (1,-1);

						\draw (-1,-2) -- (-1,1.5);
						\filldraw (-1,0) circle (2pt);
						\filldraw (-1,-1) circle (2pt);
						\filldraw (-1,1) circle (2pt);
						\draw (-1,1.5) .. controls (-0.5,1.75) .. (-0.5,2);
						
						\draw [dashed] (-1,1) -- (0,-1);
						
						\draw [loosely dotted] (-1.25,0.5) --  (-1.75,0.5);
						
						\draw (-2,-2) -- (-2,1.5);
						\filldraw (-2,0) circle (2pt);
						\filldraw (-2,-1) circle (2pt);
						\filldraw (-2,1) circle (2pt);
						\draw (-2,1.5) .. controls (-1.5,1.75) .. (-1.5,2);
						
						\draw [dashed] (-2,1) -- (-1.85,0.75); 
						\draw [dashed] (-1,-1) -- (-1.15,-0.75);
						
						\draw (-3,-2) -- (-3,2);
						\filldraw (-3,0) circle (2pt);
						\filldraw (-3,1) circle (2pt);
						\draw (-3.25,-0.75) -- (-2.75,-0.75);
						\draw (-3.25,-1) -- (-2.75,-1);
						\draw (-3.25,-1.25) -- (-2.75,-1.25);
						
						\draw [dashed] (-3,1) -- (-2,-1);
						
						\draw (-4,-2) -- (-4,2);
						\filldraw (-4,0) circle (2pt);
						\filldraw (-4,1) circle (2pt);
						\draw (-4.25,-0.75) -- (-3.75,-0.75);
						\draw (-4.25,-1) -- (-3.75,-1);
						\draw (-4.25,-1.25) -- (-3.75,-1.25);
						
						\draw [dashed] (-4,1) -- (-3,0);

						\draw [loosely dotted] (-4.25,0.5) --  (-4.75,0.5);
						
						\draw (-5,-2) -- (-5,2);
						\filldraw (-5,0) circle (2pt);
						\filldraw (-5,1) circle (2pt);
						\draw (-5.25,-0.75) -- (-4.75,-0.75);
						\draw (-5.25,-1) -- (-4.75,-1);
						\draw (-5.25,-1.25) -- (-4.75,-1.25);
						
						\draw [dashed] (-5,1) -- (-4.75,0.75);
						\draw [dashed] (-4,0) -- (-4.25,0.25);
					\end{tikzpicture}
				}
				\caption{General case}\label{fig:GStabModelRam}
			\end{subfigure}
			\hfill
			\begin{subfigure}[c]{0.45\textwidth}
				\centering	
				\resizebox{\columnwidth}{!}{%
					\begin{tikzpicture}
						\draw (6,-2) -- (6,1.5);
						\filldraw (6,0) circle (2pt);
						\filldraw (6,-1) circle (2pt);
						\filldraw (6,1) circle (2pt);
						\draw (6,1.5) .. controls (6.5,1.75) .. (6.5,2);
						
						\draw (5,-2) -- (5,1.5);
						\filldraw (5,0) circle (2pt);
						\filldraw (5,-1) circle (2pt);
						\filldraw (5,1) circle (2pt);
						\draw (5,1.5) .. controls (5.5,1.75) .. (5.5,2);
						\draw [dashed] (5,0) -- (6,0);
						\draw [dashed] (5,1) -- (6,1);

						\draw (4,-2) -- (4,1.5);
						\filldraw (4,0) circle (2pt);
						\filldraw (4,-1) circle (2pt);
						\filldraw (4,1) circle (2pt);
						\draw (4,1.5) .. controls (4.5,1.75) .. (4.5,2);
						
						\draw (3,-2) -- (3,1.5);
						\filldraw (3,0) circle (2pt);
						\filldraw (3,-1) circle (2pt);
						\filldraw (3,1) circle (2pt);
						\draw (3,1.5) .. controls (3.5,1.75) .. (3.5,2);
						\draw [dashed] (3,0) -- (4,0);
						\draw [dashed] (3,1) -- (4,1);
						
						\draw [dashed] (4,-1) -- (5,-1);  
						
						\draw (2,-2) -- (2,1.5);
						\filldraw (2,0) circle (2pt);
						\filldraw (2,-1) circle (2pt);
						\filldraw (2,1) circle (2pt);
						\draw (2,1.5) .. controls (2.5,1.75) .. (2.5,2);
						
						\draw (1,-2) -- (1,1.5);
						\filldraw (1,0) circle (2pt);
						\filldraw (1,-1) circle (2pt);
						\filldraw (1,1) circle (2pt);
						\draw (1,1.5) .. controls (1.5,1.75) .. (1.5,2);
						\draw [dashed] (1,0) -- (2,0);
						\draw [dashed] (1,1) -- (2,1);
						
						\draw [dashed] (2,-1) -- (2.25,-1);
						\draw [dashed] (3,-1) -- (2.75,-1);
						
						\draw [loosely dotted] (2.75,0.5) --  (2.25,0.5);
						
						\draw (0,-2) -- (0,2);
						\filldraw (0,0) circle (2pt);
						\filldraw (0,1) circle (2pt);
						\draw (-0.25,-0.75) -- (0.25,-0.75);
						\draw (-0.25,-1) -- (0.25,-1);
						\draw (-0.25,-1.25) -- (0.25,-1.25);

						\draw [dashed] (0,1) -- (1,-1);
						
						\draw (-1,-2) -- (-1,2);
						\filldraw (-1,0) circle (2pt);
						\filldraw (-1,1) circle (2pt);
						\draw (-1.25,-0.75) -- (-0.75,-0.75);
						\draw (-1.25,-1) -- (-0.75,-1);
						\draw (-1.25,-1.25) -- (-0.75,-1.25);
						
						\draw [dashed] (-1,1) -- (0,0);
						
						\draw (-2,-2) -- (-2,2);
						\filldraw (-2,0) circle (2pt);
						\filldraw (-2,1) circle (2pt);
						\draw (-2.25,-0.75) -- (-1.75,-0.75);
						\draw (-2.25,-1) -- (-1.75,-1);
						\draw (-2.25,-1.25) -- (-1.75,-1.25);
						
						\draw [dashed] (-2,1) -- (-1.75,0.75);
						\draw [dashed] (-1,0) -- (-1.25,0.25);

						\draw [loosely dotted] (-1.25,0.5) --  (-1.75,0.5);
						
						\draw [dashed] (-2,0) -- (-3,0) -- (-3,-3) -- (7,-3) -- (7,-1) -- (6,-1); 
						
					\end{tikzpicture}
				}
				\caption{Unramified case}\label{fig:GStabModelURam}
			\end{subfigure}
			\caption{The $\Zl$-stable curve $X^0$}		
		\end{figure}
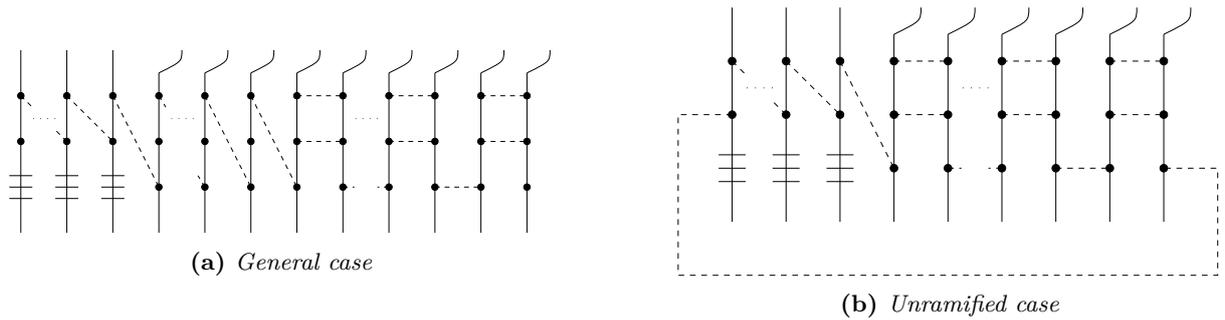

		\begin{proof}
			First, suppose $g$ is not unramified. Then by gluing along the dotted lines as in the Fig.~\ref{fig:GStabModelRam}, we obtain the desired $\Zl$-stable $C_r$-diagram $X^0$ as follows. The first part is made by gluing $p$ copies of $\Zl$-curves of genus $0$, which contributes to the $p\ell$ unramified marked points that are permuted by $\Zl$, to $1$ marked ramified point with Hurwitz data $\underline{k}=(1)$ and does not contribute to the genus. The second portion is composed of $N-2$ Seyama curves of genus $(\ell-1)/2$ glued in a chain, which contributes to $(N-2)(\ell-1)/{2}$ to the genus and to $N-2$ to the marked ramified points with free Hurwitz data. The last part is made by gluing $g'$ copies of $2$-Seyama curves. It contributes to $g'\ell$ to the genus and to $1$ ramified marked point with imposed Hurwitz data.
			
			\medskip
			
			To achieve the unramified $g$ we remove the middle section made of Seyama curves in the previous construction and glue the remaining parts on the added dotted line as in Fig.~\ref{fig:GStabModelURam}. One can easily check in the same way that it gives a desired curve.
		\end{proof}
		
		\begin{rem}
			It is readily seen that the $G$-quotient of the $G$-stable diagrams that we constructed is a $\PP^1\setminus\{0,1,\infty\}$-diagram as in \cite{iharanak}~2.1.3.  
		\end{rem}
		
		\subsection{The deformation family of $\Zl$-stable diagrams} \label{sec:deformation}
		We now start with a $\Zl$-stable $C_r$-diagram $X^0$ with $\Card \Lambda \sqcup \Lambda'\geq 2$ which is in the boundary of $\compMZlkr$ and build, by patching local formal schemes $\mW_\bullet$, $\mV_\bullet$ and $\mU_\bullet$ into a $\mS$-scheme $\mfX$ over an affine cover of $X^0$, a family of deformations $X/\Spf K[[q]]$ of $X^0$.
		
		\subsubsection{}\label{subsub:X0opens}
		Consider the following kind of families $W^0_{\lambda}$, $U^0_{\mu}$, and $V^0_{v}$ of affine open of $X^0$.
		
		\begin{enumerate}[label=\Alph*)]
			\item\label{itA:FamW} \emph{The family $(W^0_{\lambda})_{\lambda\in \Lambda}$}, resp. $(W^0_{\lambda})_{\lambda'\in \Lambda'}$, given for each $\lambda\in \Lambda$, resp. $\lambda'\in \Lambda'$, by the open complement in $X^0_{\lambda}$ of the three ramified points, resp. of the two ramified points, and represented as below:
			
			
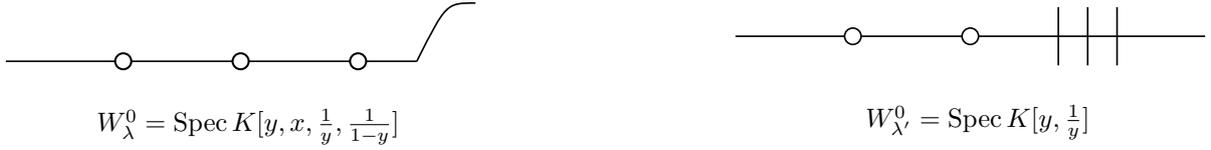
\begin{figure}[h]
				\begin{subfigure}[c]{0.4\linewidth}
					\centering	
					\resizebox{\columnwidth}{!}{%
						\begin{tikzpicture}
							\draw (-2,0) -- (1.5,0);
							\filldraw (0,0) circle (2pt);
							\filldraw [white] (0,0) circle (1.5pt);
							\filldraw (-1,0) circle (2pt);
							\filldraw [white] (-1,0) circle (1.5pt);
							\filldraw (1,0) circle (2pt);
							\filldraw [white] (1,0) circle (1.5pt);
							\draw (1.5,0) .. controls (1.75,0.5) .. (2,0.5);
						\end{tikzpicture}	
					}\\[10pt]
					{\small $W^0_{\lambda}= \Spec K[y,x, \frac{1}{y}, \frac{1}{1-y}]$}
				\end{subfigure}	
				\hfill
				\begin{subfigure}[c]{0.4\textwidth}
					\centering	
					\resizebox{\columnwidth}{!}{%
						\begin{tikzpicture}
							\draw (-2,0) -- (2,0);
							\draw[draw=black, fill=white] (0,0) circle (2pt);
							\draw[draw=black, fill=white] (-1,0) circle (2pt);
							\draw (0.75,-0.25) -- (0.75,0.25);
							\draw (1,-0.25) -- (1,0.25);
							\draw (1.25,-0.25) -- (1.25,0.25);
						\end{tikzpicture}
					}\\[10pt]
					{\small $W^0_{\lambda'}= \Spec K[y, \frac{1}{y}]$}
				\end{subfigure}
			\end{figure}
			
			\item\label{itB:FamU} \emph{The family $(U^0_{\mu})_{\mu\in M}$}, that we will specify as three subfamilies $U^0_{\mu,0,0}$, $U^0_{\mu,0,1}$, and $U^0_{\mu,1,1}$, which for $\mu\in M$ are defined such that $P^0_{\mu}$ consists of a pair of distinguished sections over $X^0_{\lambda}$ and $X^0_{\lambda'}$ with $\lambda,\lambda'\in \Lambda \sqcup \Lambda'$, and are respectively given as below (see also Fig.~\ref{Fig:ThreeSubFam}):
			
			\begin{align*}
				U^0_{\mu,0,0}&= \Spec K[y,x,y',x'][\frac{1}{y}, \frac{1}{y'}]/(T^{\lambda}_{ij} T_{ij}^{\lambda'}), &
				U^0_{\mu,0,1}&= \Spec K[x,y,y', \frac{1}{y}]/(T^{\lambda}_{ij}T^{\lambda'}_{kl}),\\
				U^0_{\mu,1,1}&= \Spec K[y,y']/(T^{\lambda}_{ij}T^{\lambda'}_{kl}).
			\end{align*}
			\begin{figure}[h]
				\begin{subfigure}[c]{0.25\textwidth}
					\centering	
					\resizebox{\columnwidth}{!}{%
						\begin{tikzpicture}
							
							\draw (-2,0) -- (1.5,0);
							\draw[draw=black, fill=white] (-1,0) circle (2pt);
							\draw[draw=black, fill=white] (1,0) circle (2pt);
							\draw (1.5,0) .. controls (1.75,0.5) .. (2,0.5);
							
							\draw (0,-2) -- (0,1.5);
							\filldraw (0,0) circle (2pt);
							\draw[draw=black, fill=white] (0,-1) circle (2pt);
							\draw[draw=black, fill=white] (0,1) circle (2pt);
							\draw (0,1.5) .. controls (0.5,1.75) .. (0.5,2);
							
						\end{tikzpicture}	
					}
					{\small $U^0_{\mu,0,0}$}
				\end{subfigure}
				\hfill
				\begin{subfigure}[c]{0.25\textwidth}
					\centering	
					\resizebox{\columnwidth}{!}{%
						\begin{tikzpicture}
							\draw (-2,0) -- (2,0);
							\draw[draw=black, fill=white] (-1,0) circle (2pt);
							\draw (0.75,-0.25) -- (0.75,0.25);
							\draw (1,-0.25) -- (1,0.25);
							\draw (1.25,-0.25) -- (1.25,0.25);
							
							\draw (0,-2) -- (0,1.5);
							\filldraw (0,0) circle (2pt);
							\draw[draw=black, fill=white] (0,-1) circle (2pt);
							\draw[draw=black, fill=white] (0,1) circle (2pt);
							\draw (0,1.5) .. controls (0.5,1.75) .. (0.5,2);
						\end{tikzpicture}	
					}
					{\small $U^0_{\mu,0,1}$}
				\end{subfigure}	
				\hfill
				\begin{subfigure}[c]{0.25\textwidth}
					\centering	
					\resizebox{\columnwidth}{!}{%
						\begin{tikzpicture}
							
							\draw (-2,0) -- (2,0);
							\draw[draw=black, fill=white] (-1,0) circle (2pt);
							\draw (0.75,-0.25) -- (0.75,0.25);
							\draw (1,-0.25) -- (1,0.25);
							\draw (1.25,-0.25) -- (1.25,0.25);
							
							\draw (0,-2) -- (0,2);
							\filldraw (0,0) circle (2pt);
							\draw[draw=black, fill=white] (0,1) circle (2pt);
							\draw (-0.25,-0.75) -- (0.25,-0.75);
							\draw (-0.25,-1) -- (0.25,-1);
							\draw (-0.25,-1.25) -- (0.25,-1.25);
							
						\end{tikzpicture}	
					}
					{\small $U^0_{\mu,1,1}$}
				\end{subfigure}
				\caption{The three subfamilies of $U^0_{\mu}$}\label{Fig:ThreeSubFam}
			\end{figure}
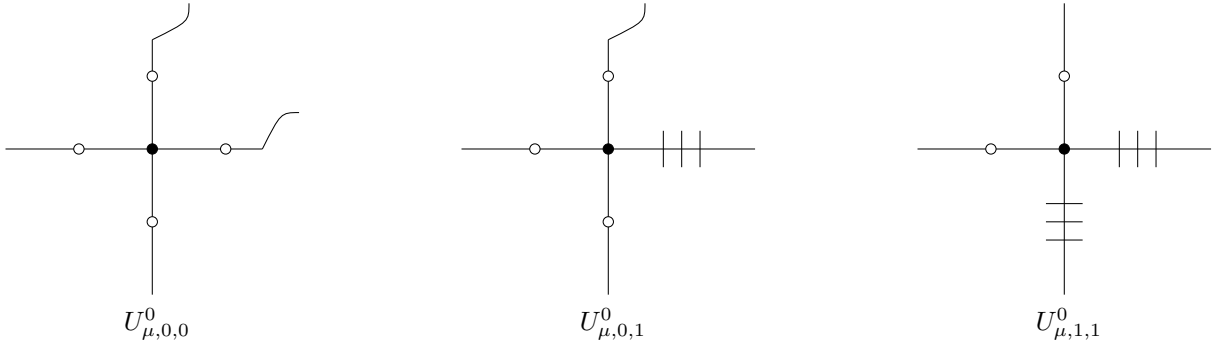
			
			\item\label{itC:FamV} \emph{The family $(V^0_{v})_{v\in N}$} given for each $v\in N$ by taking the component $X^0_{\lambda}$ that supports the section $Q^0_v$ and removing all the other distinguished sections, to obtain
			\begin{equation}
				V^0_v= \Spec K[y,x,\frac{1}{1-y}],\text{ resp. }V^0_v= \Spec K[y].
			\end{equation}
			for $\lambda\in \Lambda$, resp. $\lambda\in \Lambda'$.
		\end{enumerate}
		
		\medskip
		
		We thus obtain an affine cover of the $\Zl$-stable curve $X^0$
		\[
		X^0= \bigcup_{\lambda\in \Lambda\sqcup \Lambda'} W^0_{\lambda} \bigcup_{\mu \in M} U^0_{\mu} \bigcup_{v\in N} V^0_v
		\]
		where each open is $\Zl$-stable by construction, and such that:
		
		\begin{enumerate}
			\item For $\mu \in M$ such that $P^0_{\mu}$ contains a distinguished section of $X^0_{\lambda}$ and $X^0_{\lambda'}$ we have $W^0_{\lambda}$ and $W^0_{\lambda'}$ as open subsets of $U^0_{\mu}$ and $W^0_{\lambda}\cap W^0_{\lambda'}=\varnothing$.
			
			\item For $v\in N$ such that $Q^0_v$ is a distinguished section of $X^0_{\lambda}$ we have $W^0_{\lambda}$ as an open subset of $V^0_v$. 
			
			\item The intersection of $U^0_{\mu}$ or $V^0_v$ with any other member of the affine cover is either empty, $W^0_{\lambda}$ or $W^0_{\lambda}\sqcup W^0_{\lambda'}$.
		\end{enumerate}	
		
		These properties ensure, in the next section, the possibility of patching local formal schemes over the affine cover that we just defined.
		
		\subsubsection{}
		Consider the affine formal scheme $\mS= \Spf K[[q]]$ with ideal of definition $\mfq=(q)$ equipped with the $G$-action $q\mapsto \zeta_{\ell}q$ by our choice of isomorphism $G\simeq \mu_{\ell}$. In order to construct a formal scheme $\mfX$ with base $X^0$ over $\mS$ with a compatible $G$-action we shall define affine formal $\mS$-schemes $\mW_{\lambda}$, $\mU_{\mu}$ and $\mV_v$ with bases $W^0_{\lambda}$, $U^0_{\mu}$ and $V^0_v$ with ideal of definition the pullbacks of $\mfq$ denoted $\mfq$ again.
		
		\medskip
		
		For $\lambda,v$ we set 
		\begin{equation*}
			\mW_{\lambda}= \Spf \Gamma(W_{\lambda}^0,\OO_{X^0})[[q]]\text{ and }\mV_v=\Spf \Gamma(V_{v}^0,\OO_{X^0})[[q]],
		\end{equation*}
		where the $K$-algebras of sections $\Gamma(-,\OO_{X^0})$ are given by one of the explicit $K$-algebra of the affine schemes of \ref{subsub:X0opens} \ref{itA:FamW} and \ref{itC:FamV} above.

		\medskip
		
		Whenever $\lambda$ and $v$ are such that $W^0_{\lambda}$ is an open subset of $V^0_v$ the open immersion $j_{v/\lambda}\colon \mW_{\lambda}\rightarrow \mV_v$ over it is obtained without effort. For instance, let us assume $W^0_{\lambda}=\Spec K[y,x,\frac{1}{y},\frac{1}{1-y}]$ and $V^0_v= \Spec K[y,x,\frac{1}{1-y}]$. By \cite{ega1} Proposition~10.1.4 it suffices to check that the map $K[y,x,\frac{1}{1-y}][\frac{1}{y}][[q]]^{\wedge} \rightarrow \Gamma(\mW_{\lambda})$, where $\wedge$ denotes the $\mfq$-adic completion, is an isomorphism. But this is clear by construction. Note that $j_{v/\lambda}$ is an $\mS$-morphism. 
		
		\medskip
		
		Let us consider the case of $\mU_\mu$, whose base $U^0_\mu$ is obtained via 3 subfamilies $U^0_{\mu,0,0}$, $U^0_{\mu,0,1}$, and $U^0_{\mu,1,1}$ as in Section~\ref{subsub:X0opens} \ref{itB:FamU}.
		\begin{prop} \label{openjmu}
			For $\mu\in M$ such that $U^0_{\mu}$ is of the form $U^0_{\mu,0,0}$. Let us define
			\[\textstyle
			\mU_{\mu,0,0}=\Spf K[T,T',X,X'][\frac{1}{1-T},\frac{1}{1-T'}][[q]]/(T_{r}T'_{r'}-q)\text{ with }\begin{cases}
				\textstyle X^{\ell}=T^r(1-T), &\textstyle T_{r}=\zeta_{2\ell}X\\
				\textstyle {X'}^{\ell}={T'}^{r'}(1-T'), &\textstyle T'_{r'}=\zeta_{2\ell}{X'}.\\
			\end{cases}
			\]
			Then we can identify $\mU_{\mu,0,0} \mod \mfq$ with $U^0_{\mu,0,0}$ by $(T,T')\mapsto (y,y')$ with the choices $(X,X')\mapsto (x,x')$. Furthermore, for $\lambda \in \Lambda$ such that $W^0_{\lambda}=\Spec K[y,x,\frac{1}{y},\frac{1}{1-y}]$, the scheme $W^0_{\lambda}$ is an open subset of $U^0_{\mu,0,0}$ given by inverting $y$, so that $\Gamma(\mU_{\mu})[\frac{1}{T}]^{\wedge}\rightarrow \Gamma(\mW_{\lambda})$, given by $T\rightarrow y$, is an isomorphism, which induces an open immersion $j_{\mu/\lambda}\colon \mW_{\lambda} \rightarrow \mU_{\mu,0,0}.$ 
		\end{prop}
		
		\begin{proof}
			By assumption $T_{r}T'_{r'}=q$ so that for $N\geq 1$ we have
			\[
			\Gamma(\mU_{\mu})[\frac{1}{T}]/\mfq^N= K[T,T',X,X'][\frac{1}{1-T},\frac{1}{1-T},\frac{1}{T'}]/(T_{r}T'_{r}-q)^N.
			\]	
			
			As $T$ and $1-T$ are invertible, so is $X^{\ell}=T^r(T-1)$, and thus so is $X$ and $T_{r}=\zeta_{2\ell}X$. It follows that $(T_{r}T'_{r'})^N=0$ if and only if ${T'_{r'}}^N=0$. Now as ${T'_{r'}}^{\ell N}= -{T'}^{r'N} (1-T')^{N}$ we have ${(1-T')}^N=0$, which gives that ${T'}^{-1}$ can be written as $\sum\nolimits_{k=0}^{N-1} {(1-T')}^k$. To recover $T'$ and $X'$ first as ${T'}^{\ell}_{r'}=-{T'}^{r'}(T'-1)$ we have ${T'}^{\ell}_{r'}=P(T')$ with $P$ invertible for the composition in $K[[T'-1]]$. So there is $F\in K[[T'-1]]$ such that $F({T'}^{\ell}_{r'})=T'$. As ${T'}^{\ell}_{r'}$ is nilpotent of order $N$ we can truncate $F$ to get a polynomial $\widetilde{F}$ that verifies the equality $\widetilde{F}({T'}^{\ell}_{r'})=T'$ in $\Gamma(\mU_{\mu})[\frac{1}{T}]/\mfq^N$. 
			Thus, as $T'_{r'}=\frac{q}{\zeta_{2\ell} T}$, we have 
			\[
			\Gamma(\mU_{\mu})[\frac{1}{T}]/\mfq^N=K[X,T,\frac{1}{T},\frac{1}{1-T}][q]/(q^N)
			\]
			and the desired isomorphism by passing to the $\mfq$-adic completion. It is clear that this isomorphism is compatible with the $\Zl$-action on both sides.
		\end{proof}
		
		The other open immersions are proven in the same way. 
		
		\subsubsection{} \label{subsub:formaldef}
		One thus obtains a proper formal regular $\mS$-scheme $\mfX$ with a collection of sections $(\mQ_v)_{v\in N}$ with base space $X^0$ by gluing along the affine formal schemes $\mU_{\mu},\mV_v$ and $\mW_{\lambda}$. 
		
		\medskip
		
		The formal scheme $\mfX$ has the property that, for each $\mu$, $\lambda$ or $v$ we have $\mS$-isomorphisms
		\[
		\varphi_{\mu}\colon\mfX_{|U^0_{\mu}}\rightarrow \mU_{\mu},~\varphi_{\lambda}\colon \mfX_{|W^0_{\lambda}}\rightarrow \mW_{\lambda}, \varphi_v\colon \mfX_{|V^0_v}\rightarrow \mV_v
		\]
		extending the identity maps of $U^0_{\mu}$, $W^0_{\lambda}$ and $V^0_v$, respectively, such that
		\begin{enumerate}[label=(\alph*)]
			\item for each $v\in N$, $\mQ_v$ is induced from the canonical section $\mS\rightarrow \mfX_{|V^0_v}$ that lift the section $Q_v$ of $V^0_v$,
			
			\item the isomorphisms $\varphi_{\lambda}$, $\varphi_{\lambda}$ and $\varphi_v$ are compatible with the open immersions $j_{v/\lambda}$ and $j_{\mu/\lambda}$. 
		\end{enumerate}

		\medskip
		
		A direct application of Grothendieck's existence theorem \cite{EGA3.1}~5.4 as in \cite{iharanak}~2.4.1 and 3.1 provides the algebraization of the formal scheme $\mfX$ into a scheme $X$ over $\Spec K[[q]]$, whose generic fiber $X_{\eta}$ is a smooth geometrically irreducible genus $g$ curve with $m$ marked points and a $\Zl$-action, coming by pullback of the one on $X$, with Hurwitz data $\underline{kr}$, and whose special fiber is $X^0$.
		
		\medskip
		
		The sets of divisors $D=((X^0_{\lambda})_{\lambda\in \Lambda\sqcup \Lambda'}, (Q_v)_{v\in N})$ and $\mD=((X^0_{\lambda})_{\lambda\in \Lambda\sqcup \Lambda'}, (\mQ_v)_{v\in N})$ are regular with normal crossings on $X$ and $\mfX$ respectively in the sense of \cite{gm} Section~1.8.3, see \cite{iharanak}~3.2 for details. 
		
		\begin{rem}
			The generic fiber $X_{\eta}$ of the scheme $X$ should be interpreted as a tangential base point $\eta\colon \Spec K((q))\rightarrow \MZlkr$ in the moduli space.
		\end{rem}
		
		\subsubsection{}\label{subsub:DanDprime}
		Another important output of our construction, that will be of interest in the next section, is that we can explicitly track our tangential base points in the different formal completions of $\mfX$ along chosen closed subsets of the special fiber. 
		
		\medskip
		
		Consider the completion $\mfX_{\mu}$ of $\mfX$ along $P^0_{\mu}$. By construction $\mfX_{| U^0_{\mu}}= \Spf A/(T_{r,s}T'_{r',s'} -q)$ for a ring $A$ given in the construction of $U^0_{\mu}$ and $P^0_{\mu}$ corresponds to the ideal $(T_{r,s},T'_{r',s'})$, so that
		\[
		\mfX_{\mu} = \Spf K[[T_{ij},T'_{kl}]] \text{ with } T_{ij}T'_{kl}=q \text{ as usual}. 
		\]
		Let $T_1,T_2$ be two indeterminacies. We have a commutative diagram
		\[
		\begin{tikzcd}
			\Spf K[[T_1,T_2]] \arrow[d, "q\mapsto T_1T_2"'] \arrow[r, "\mu"] &  \mfX_{\mu} \arrow[d] \\
			\Spf K[[q]] \arrow[r, "s"]  & \mS
		\end{tikzcd}    
		\]
		where $K[[T_1,T_2]]$ has ideal of definition $(T_1T_2)$, and where the top horizontal map is an isomorphism. The formal scheme $\mfX_{\mu}$ comes with a divisor $\mD_{\mu}$ given by the pullback of $\mD$ which has two components corresponding to $X^0_{\lambda}$ and $X^0_{\lambda'}$ where $\lambda(\mu)=(\lambda,\lambda')$. They are defined by $T_1=0$ and $T_2=0$ respectively, so that $\mD_{\mu}$ is a set of divisors with regular normal crossing on $\mfX_{\mu}$. 
		
		\medskip
		
		We shall also consider the completion $\mfX_{\lambda}$ of $\mfX$ along $X^0_{\lambda}$. It is also equipped with a divisor $\mD_{\lambda}$ as the pullback of $\mD$ to $\mfX_{\lambda}$ which consists of the union of two divisors:
		\begin{enumerate}
			\item $\mD^0_{\lambda}$ given by $X^0_{\lambda}$
			\item $\mD'_{\lambda}$ given by the distinguished sections of $X^0_{\lambda}$. 
		\end{enumerate}
		It is again a set of divisors with regular normal crossings. 
		
		\medskip
		
		By arguing as in the proof of Proposition~\ref{openjmu}, one further obtain the following compatibility result between tangential base points and formal completions.
		
		\begin{prop}\label{pro:formaldiagram}
			Let $\mu\in M$ and $\lambda\in \lambda(\mu)$. Then we have the following commutative diagram in the category of formal schemes
			$$\begin{tikzcd}
				\mfX_{\mu} \arrow[dd] \arrow[dddr] &  \Spf K[[T_1,T_2]] \arrow[l, "\mu"'] \arrow[ddd]  &\Spf K((T))[[q]]   \arrow[l]  & \Spf K((T)) \arrow[dd, "T^{\lambda}_{ij}"] \arrow[l] \arrow[dl] \\
				& & \Spf K[[T]] \arrow[ul] \arrow[ddr, crossing over] & \\
				\mU_{\mu} \arrow[d, hook] &  & & W^0_{\lambda} \arrow[d, hook] \arrow[lll, hook, crossing over] \\
				\mfX & \mfX_{\lambda} \arrow[l] &  & X^0_{\lambda} \arrow[ll, hook]
			\end{tikzcd}$$
		\end{prop}
		
		One remarks that the map $\iota_{\lambda}\colon \Spec K[[T]] \rightarrow \Spf K[[T_1,T_2]]$ which is given by the quotient by $T_2$ factors through the restriction to the special fiber $\Spec K[[T_1,T_2]]/(T_1T_2)$.
		
		\section{Galois actions by Grothendieck-Murre theory} \label{sec:galoisact}
		
		Starting with a $G$-stable diagram $X^0$ with $\Card \Lambda\sqcup \Lambda'\geq 2$ the end result of the previous section gives us a smooth curve $X_{\eta}$ represented by a $K((q))$-point on $\mM_{g,[m]}(G)_{\underline{kr}}^\nu$ which comes with a model $X$ over $\mS$ with special fiber $X^0$. 
		
		\medskip
		
		We will now relate the Galois action on the fundamental groupoid $\Pi_1(X_{\eta} \setminus \{ (Q_v)_{v\in N}\}, (\vec{\mu})_{\mu\in M})$ of $X_{\eta}$ based at the punctures coming from the double points $(P_{\mu})_{\mu\in M}$ of $X$ to the ones on the curves $(C_r)_{r\in \{0,\dots, \ell-2\}}$ obtained by the tangential base points $T^r_{ij}$ that we defined in Section \ref{sec:seyamacurves}. To do so, we follow some equivalence between categories of covers as in \cite{iharanak}: the category $\operatorname{Rev}^{D}(X)$ of finite étale covers of $X$ tamely ramified along the divisor $D$, made of the union $X^0\cup \{(Q_v)_{v\in N}\}$, is canonically equivalent both to $\operatorname{Rev}^{\mD}(\mfX)$ and $\operatorname{Rev}(X_{\eta} \setminus \{ (Q_v)_{v\in N}\})$. For $\mu\in M$, we then define some fiber functors $\vec{\mu}$, so that, by the previous canonical equivalences of categories, we have the isomorphism
		\begin{equation*}
			\pi_1^{\mD}(\mfX, \vec{\mu}) \simeq \pi_1^{D}(X, \vec{\mu}) \simeq \pi_1(X_{\eta} \setminus \{(Q_v)_{v\in N}\}, \vec{\mu}).
		\end{equation*} 
		
		\medskip
		
		Those equivalences are Galois equivariant so in order to determine whenever an element of $G_K$ acts trivially on the geometric part of $\pi_1(X_{\eta} \setminus \{(Q_v)_{v\in N}\}, \vec{\mu})$ it is enough to do so on the left-hand side.
		
		\subsection{Tamely ramified fundamental groups and fiber functors}
		
		\subsubsection{}
		We start by defining fiber functors on $\Rev^{\mD}(\mfX)$ locally by fixing $\mu\in M$ and considering $\mfX_{\mu}$. Recall that we have a commutative diagram
		\[\begin{tikzcd}
			\Spf K[[T_1,T_2]] \arrow[r, "\mu"] \arrow[d] & \mfX_{\mu} \arrow[d] \\[-5pt]
			\Spf K[[q]]\arrow[r, "s"'] & \mS
		\end{tikzcd}\]
		given by the map $q\mapsto T_1T_2$. Both maps $s$ and $\mu$ define fiber functors, $\vec{\mu}$ for $\Rev^{\mD_{\mu}}(\mfX_{\mu}) $ and $\vec{s}$ for $\Rev^{S^0}(\mS)$, see \cite{iharanak} 3.3.1 and 3.3.2. 
		
		\medskip
		
		To be explicit, consider a compatible choice of indeterminates $\{T_1^{\frac{1}{N}},T_2^{\frac{1}{N}}\}_{N\in \N}$ and $\{q^{\frac{1}{N}}\}_{N\in \N}$ to form the fields $K\{\{T_1,T_2\}\}$ and $K\{\{q\}\}$. Then for $\mfB=\Spf \mB \in \Rev^{\mD}(\mfX_{\mu})$, resp. $\mfA=\Spf \mA \in \Rev^{S^0}(\mS)$, the value of the fiber functors are given by
		\begin{equation*}
			\vec{\mu}(\mfB)=\Hom_{K[[T_1,T_2]]}(\mB, \overline{K}\{\{T_1,T_2\}\}),\text{ resp. } \vec{s}(\mfA)= \Hom_{K[[q]]}(\mA, \overline{K}\{\{q\}\}).
		\end{equation*}
		
		\subsubsection{}
		By choosing geometric points such that $q^{\frac{1}{N}} \mapsto (T_1T_2)^{\frac{1}{N}}$, one obtains two compatible homotopy exact sequences 
		\begin{equation}\label{eq:homhotTwice}
			\begin{tikzcd}
				1 \arrow[r] & \widehat{\Z}(1)\times \widehat{\Z}(1) \arrow[d] \arrow[r, "j_{\mu}"] & \pi_1^{\mD_{\mu}}(\mfX_{\mu}, \vec{\mu}) \arrow[r, "p_{\mu}"] \arrow[d, "p_{\mu/\mS}", swap] & G_K \arrow[r] \arrow[d, equal] \arrow[l, "s_{\mu}", bend right = 50] & 1 \\
				1 \arrow[r] &  \widehat{\Z}(1)  \arrow[u, shift left=2.5, "j_{\lambda}"] \arrow[u, shift right=2.5, "j_{\lambda'}", swap] \arrow[r, "j_{\mS}"] & \pi_1^{S^0}(\mS, \vec{s}) \arrow[r, "p_{\mS}"] & G_K \arrow[r] \arrow[l, "s_{s}", bend right = 50] & 1
			\end{tikzcd}
		\end{equation}
		where the geometric parts $\widehat{\Z}(1)\times \widehat{\Z}(1)$ and $\widehat{\Z}(1)$ are equipped with the Galois actions coming from the sections defined by the choices of tangential base points $\mu$ and $s$. We refer to \cite{iharanak} 3.3.1-3.3.4 for details.
		
		\subsubsection{}
		We will now track explicitly the fiber functors defined by $\vec{\mu}$ on $\Rev^{\mD_{\lambda}}(\mfX_{\lambda})$ and $\Rev^{\mD_{\lambda'}}(\mfX_{\lambda'})$ for $(\lambda, \lambda')=\lambda(P^0_{\mu})$ and compare them to the one given by the tangential base points $T^{\lambda}_{ij}$ of Definition~\ref{def:Gdiagram}. First of all, remark that the map $\mfX_{\mu} \rightarrow \mfX_{\lambda}$ pulls back the divisor $\mD_{\lambda}$ to $\mD_{\mu}$ so that it induces a base change functor
		\[
		\Rev^{\mD_{\lambda}}(\mfX_{\lambda}) \longrightarrow \Rev^{\mD_{\mu}}(\mfX_{\mu}),
		\]
		and thus we have a fiber functor on $\Rev^{\mD_{\lambda}}(\mfX_{\lambda})$ that is given by composition with $\vec{\mu}$, which we also denote by~$\vec{\mu}$. In particular this comes with a map on the étale fundamental groups
		\[
		p_{\mu/\lambda}\colon \pi_1^{\mD_{\mu}}(\mfX_{\mu}, \vec{\mu})\longrightarrow \pi_1^{\mD_{\lambda}}(\mfX_{\lambda}, \vec{\mu}).
		\]
		
		In the same way, the morphism $f_{\lambda}\colon \mfX_{\lambda}\rightarrow \mS$ defines a map
		\[
		p_{\lambda/\mS}\colon \pi_1^{\mD_{\lambda}}(\mfX_{\lambda}, \vec{\mu}) \longrightarrow \pi_1^{S^0}(\mS, \vec{s})
		\]
		by the fact that the pullback of $S^0$ is the divisor $\mD_{\lambda}^0\cup \mD_{\lambda}''$ where $\mD_{\lambda}''$ is given by $\mD_{\lambda}'$ restricted to~$X^0_{\lambda}$. As the map $\mfX_{\mu} \rightarrow \mfX_{\lambda}$ is a map of $\mS$-schemes we have the commutativity condition
		\[
		p_{\lambda/\mS} \circ p_{\mu /\lambda} = p_{\mu/ \mS}
		\]	
		and compatibility with the previous homotopy exact sequences of Eq.~\eqref{eq:homhotTwice}.
		
		\subsubsection{}\label{subsub:equiCatRevFormSpec}
		By Theorem~4.3.2 of \cite{gm} the restriction map to $X^0_{\lambda}$ gives a categorical equivalence
		\[
		\Rev^{\mD_{\lambda}'}(\mfX_{\lambda}) \simeq \Rev^{D_{\lambda}}(X^0_{\lambda}),
		\]
		and the last one is canonically equivalent to $\Rev( W_0^{\lambda})$.

		\begin{prop}
			The isomorphisms $\Rev^{\mD'_{\lambda}}(\mfX_{\lambda}) \simeq \Rev^{\mD_{\lambda}}(X^0_{\lambda}) \simeq \Rev (W^0_{\lambda})$ transform the fiber functor $\vec{\mu}$ in $\overrightarrow{T_{ij}}$ and thus yields a Galois compatible isomorphism
			\[
			\pi_1^{\mD'_{\lambda}}(\mfX_{\lambda}, \vec{\mu}) \simeq \pi_1(W^0_{\lambda}, \overrightarrow{T}_{ij}).
			\]
		\end{prop}
		
		\begin{proof}
			By Proposition~\ref{pro:formaldiagram}, the following diagram commutes
			\[\begin{tikzcd}
				\Spf K[[T_1,T_2]] \arrow[d, "\mu"'] &  & \Spf K[[T]] \arrow[ll, "\iota_{\lambda}"'] \arrow[d, "T_{ij}"] \\
				\mfX_{\mu} \arrow[r] & \mfX_{\lambda} & X^0_{\lambda} \arrow[l]
			\end{tikzcd}\]
			where we recall the map $\iota_{\lambda}\colon \Spec K[[T]] \rightarrow \Spf K[[T_1,T_2]]$ is given by the quotient by $T_2$.
			
			It thus suffices to check that the fiber functors on $\Rev^{\mD_{\lambda}'}(\mfX_{\lambda})$ given by $\mu$ and $\mu\circ \iota_{\lambda}$ are canonically equivalent and that they are also equivalent to the one given by composition of the pullback to the special fiber and  $\overrightarrow{T}_{ij}$.
			
			Let $\mfB\in \Rev^{\mD_{\lambda}'}(\mfX_{\lambda})$ and consider $A\in \Rev^{D_{\lambda}}(X^0_{\lambda})$ obtained from $\mfB$ by base change to the special fiber. The pullback of $\mfB$ to $\mfX_{\mu}$ is $\Spf \mB \in \Rev^{(T_1=0)}( \mfX_{\mu})$ with $\mB$ a direct sum of subalgebras of $\overline{K}[[T_1^{\frac{1}{N}}, T_2]]$ for some $N\geq 1$. 
			Then we have
			\begin{align*}
				\vec{\mu}(\mfB)&= \Hom_{K[[T_1,T_2]]}(\mB, \overline{K}\{\{T_1,T_2\}\})\\
				&=\Hom_{K[[T_1,T_2]]} ( \mB, \overline{K}\{\{T_1\}\} [[T_2]]) \\
				&= \Hom_{K[[T]]}(\mB/T_2, \overline{K}\{\{T\}\}) \\
				\vec{\mu}(\mfB)&= \overrightarrow{\mu\circ \iota_{\lambda}}(\mfB) =\overrightarrow{T_{ij}}(A).
			\end{align*} 
			
		\end{proof}
		
		\begin{rem}
			The map $\Spf K[[T]] \rightarrow \Spf K[[T_1,T_2]]$ does not define a base change $\Rev^{\mD_{\mu}}(\mfX_{\mu})\rightarrow \Rev^{(T=0)}( \Spf K[[T]])$ as the pullback of the divisor $\mD_{\mu}$ is $\Spf K[[T]]$ and not $(T=0)$. Thus we can not define a fiber functor for the first category in this way. 
		\end{rem}
		
		\subsubsection{}
		We can now compare $\pi_1^{\mD_{\lambda}'}(\mfX_{\lambda}, \vec{\mu})$ and $\pi_1^{\mD_{\lambda}}(\mfX_{\lambda}, \vec{\mu})$ by Grothendieck-Murre theory since $\mD_{\lambda}$ and $\mD_{\lambda}'$, as defined in Section~\ref{subsub:DanDprime}, are two divisors that differ by the special fiber, see \cite{gm}~Corollary~5.1.11.
		
		\begin{prop} \label{pro:exactxlambda}
			We have an exact sequence
			\[\begin{tikzcd}
				1\arrow[r] & \widehat{\Z}(1)\arrow[r, "\alpha"] & \pi_1^{\mD_{\lambda}}(\mfX_{\lambda}, \vec{\mu}) \arrow[r] & \pi_1^{\mD_{\lambda}'}(\mfX_{\lambda}, \vec{\mu}) \arrow[r] & 1
			\end{tikzcd}\]
			where $\alpha= p_{\mu/\lambda}\circ j_{\mu}\circ j_{\lambda}$ and where $\beta$ comes from the canonical projection induced by the inclusion $\Rev^{\mD_{\lambda}'}(\mfX_{\lambda})\subset \Rev^{\mD_{\lambda}}(\mfX_{\lambda})$.  
		\end{prop}
		
		\begin{proof}
			By \cite{gm} Theorem 7.3.1 we have the exactness of the sequence 
			\[\begin{tikzcd}
				\widehat{\Z}(1)\arrow[r, "\alpha"] & \pi_1^{\mD_{\lambda}}(\mfX_{\lambda}, \vec{\mu}) \arrow[r] & \pi_1^{\mD_{\lambda}'}(\mfX_{\lambda}, \vec{\mu}) \arrow[r] & 1.
			\end{tikzcd}
			\]
			The injectivity of $\alpha$ can be deduced from the injectivity of $p_{\lambda/\mS} \circ \alpha= j_{\mS}$. 
		\end{proof}
		
		\begin{rem}
			With the equality $p_{\lambda/\mS}\circ p_{\mu /\lambda} \circ s_{\mu} = s_{s}$ we also have the surjectivity of $p_{\lambda/\mS}$. 
		\end{rem}

		\subsection{Geometric Galois actions and groupoids}
		
		For the fundamental group of a curve $X$ over $K$ the geometric part is defined to be the fundamental group of $X_{\overline{K}}$ and coincide with the kernel of the projection to $G_K$ given by the arithmetic geometric fundamental homotopy exact sequence Eq.~\eqref{eq:FAGexseq}.
		
		\subsubsection{} Following \cite{iharanak} 3.4.7 we define  geometric parts of the fundamental groups $\pi_1^{\mD_{\lambda}} (\mfX_{\lambda}, \vec{\mu})$ as the kernels of such projections to $G_K$.
		
		\begin{defin} \label{def:geopart}
			The geometric part $\pi_1^{\mD_{\lambda}}({\mfX_{\lambda}}_{\overline{K}}, \vec{\mu})$ of $\pi_1^{\mD_{\lambda}} (\mfX_{\lambda}, \vec{\mu})$ is the kernel of $p_{\lambda}=p_{\mS}\circ p_{\lambda/\mS}$.
		\end{defin}
		
		\begin{prop}\mbox{} \label{pro:structfondgroup} We have the following results on the structure of $\pi_1^{\mD_{\lambda}}(\mfX_{\lambda}, \vec{\mu})$.
			\begin{enumerate}
				\item \label{item:GMGensplit} We have an exact sequence
				\begin{equation*}
					\begin{tikzcd}
						1\arrow[r] &\widehat{\Z}(1)\arrow[r, "\overline{\alpha}"] & \pi_1^{\mD_{\lambda}}({\mfX_{\lambda}}_{\overline{K}}, \vec{\mu}) \arrow[r] & \pi_1^{D_{\lambda}}({X^0_{\lambda}}_{\overline{K}}, \overrightarrow{T_{ij}}) \arrow[r] & 1
					\end{tikzcd}
				\end{equation*}
				and an isomorphism $\pi_1^{\mD_{\lambda}}({\mfX_{\lambda}}_{\overline{K}}, \vec{\mu}) \simeq \widehat{\Z}(1)\times \pi_1^{D_{\lambda}}({X^0_{\lambda}}_{\overline{K}}, \overrightarrow{T_{ij}})$.
				
				\item The exact sequence 
				\[\begin{tikzcd} 
					1 \arrow[r] & \pi_1^{\mD_{\lambda}}({\mfX_{\lambda}}_{\overline{K}}, \vec{\mu}) \arrow[r] & \pi_1^{\mD_{\lambda}}(\mfX_{\lambda}, \vec{\mu}) \arrow[r] & G_K\arrow[r] & 1
				\end{tikzcd}
				\]
				admits a splitting and we have an isomorphism $\pi_1^{\mD_{\lambda}}(\mfX_{\lambda}, \vec{\mu})\simeq \pi_1^{\mD_{\lambda}}({\mfX_{\lambda}}_{\overline{K}}, \vec{\mu}) \rtimes G_K$. 
			\end{enumerate}
			Furthermore, the action of $G_K$ on $\pi_1^{\mD_{\lambda}}({\mfX_{\lambda}}_{\overline{K}}, \vec{\mu})$ preserves the direct product decomposition of \ref{item:GMGensplit} and induces the Galois action on $\pi_1^{D_{\lambda}}({{X^0_{\lambda}}}_{\overline{K}}, \overrightarrow{T_{ij}})$ given by the tangential base point $T^{\lambda}_{ij}$. 
		\end{prop}
		
		\begin{proof}\mbox{} 
			\begin{enumerate}
				\item We deduce the exact sequence from the one of Proposition~\ref{pro:exactxlambda}, where we replaced the last term via the equivalence of categories $\Rev^{\mD_{\lambda}'}(\mfX_{\lambda}) \simeq \Rev^{D_{\lambda}}(X^0_{\lambda})$, see \ref{subsub:equiCatRevFormSpec}.
				
				We know that $\widehat{\Z}(1)$ is the kernel of $\beta$ so that its image lands in the geometric part is a given. The short exact sequence follows.
				
				The projection $p_{\lambda/\mS}$ induces a geometric counterpart 
				$$\overline{p_{\lambda/\mS}}\colon \pi_1^{\mD_{\lambda}}({\mfX_{\lambda}}_{\overline{K}}, \vec{\mu}) \longrightarrow \widehat{\Z}(1)$$
				which verifies $\overline{p_{\lambda/\mS}} \circ \overline{\alpha} =\mathrm{id}_{\widehat{\Z}(1)}$. It follows that $\Ker \overline{p_{\lambda/\mS}} \cap \overline{\alpha}(\widehat{\Z}(1))=\{1\}$ so that $\Ker \overline{p_{\lambda/\mS}}$ is isomorphic to $\pi_1^{D_{\lambda}}({X^0_{\lambda}}_{\overline{K}}, \overrightarrow{T_{ij}})$ and we have the direct product decomposition. 
				
				\item The splitting is given by $s_{\mu}\circ p_{\mu/\lambda}$. The fact that the resulting $G_K$-action preserves the direct product decomposition and induces the $G_K$-action on $\pi_1^{D_{\lambda}}({X^0_{\lambda}}_{\overline{K}}, \overrightarrow{T_{ij}})$ given by the tangential base point $T_{ij}$, follows directly from the compatibility of the fiber functors $\vec{\mu}$, $\vec{s}$ and $\overrightarrow{T_{ij}}$.
			\end{enumerate}
		\end{proof}
		
		\subsubsection{} We can now state the basic result that determines when an element of $G_K$ acts trivially on $\pi_1^{\mD_{\lambda}}({\mfX_{\lambda}}_{\overline{K}}, \vec{\mu}).$

		\begin{prop}\label{pro:trivialactiongroup}
			An element of $G_K$ acts trivially on $\pi_1^{\ell,\mD_{\lambda}}({\mfX_{\lambda}}_{\overline{K}}, \vec{\mu})$ if and only if it acts trivially on $\pi_1^{\ell,D_{\lambda}}({X^0_{\lambda}}_{\overline{K}}, \overrightarrow{T_{ij}}).$
		\end{prop}
		\begin{proof}
			The decomposition of $\pi_1^{\mD_{\lambda}}({\mfX_{\lambda}}_{\overline{K}}, \vec{\mu})$ given by \ref{item:GMGensplit} of the previous result passes to the pro-$\ell$-completion, which gives
			\[
			\pi_1^{\ell,\mD_{\lambda}}({\mfX_{\lambda}}_{\overline{K}}, \vec{\mu})\simeq \widehat{\Z_{\ell}}(1) \times \pi_1^{\ell,D_{\lambda}}({X^0_{\lambda}}_{\overline{K}}, \overrightarrow{T_{ij}}).
			\]
			As the $G_K$-action preserves the product, the implication is straightforward. For the reciprocal, let $\sigma \in G_K$ that acts trivially on $\pi_1^{\ell,\mD_{\lambda}}({\mfX_{\lambda}}_{\overline{K}}, \vec{\mu})$. Let us choose a representation $(y_1,\dots, y_{2g}, x_1,\dots x_n \mid \prod\nolimits_i [y_i,y_{i+1}] x_1\cdots x_{n} )$ of $\pi_1^{\ell,D_{\lambda}}({X^0_{\lambda}}_{\overline{K}}, \overrightarrow{T_{ij}})$ in the usual way, where $x_1$ denotes the loop around the closed point image of $T_{ij}$ in $X^0_{\lambda}$. We have $\sigma(x_1)=x_1^{\chi_{\ell}(\sigma)}=x_1$ by assumption. But $\sigma$ also acts by $\chi_{\ell}(\sigma)$ on the first factor $\widehat{\Z_{\ell}}(1)$ so the action of $\sigma$ on $\pi_1^{\ell,\mD_{\lambda}}({\mfX_{\lambda}}_{\overline{K}}, \vec{\mu})$ is trivial. 
		\end{proof}
		
		\begin{rem}
			More generally, the result also holds in the case of any almost full class of finite groups $\mC$ and the maximal pro-$\mC$-quotients of $\pi_1^{\mD_{\lambda}}({\mfX_{\lambda}}_{\overline{K}}, \vec{\mu})$ and $\pi_1^{D_{\lambda}}({X^0_{\lambda}}_{\overline{K}}, \overrightarrow{T_{ij}})$, see \cite{iharanak} Proposition~3.4.8. 
		\end{rem}
		
		\subsubsection{} In order to conclude, we first we need to explain how to move from fundamental groups to fundamental groupoids. This is essentially formal and comes down to the fact that the set of étale paths are principal homogeneous spaces under the translation actions of the fundamental groups. As such, the technical details will mostly be avoided.
		
		\medskip
		
		Let $M_{\lambda}=\{\mu \in M \mid \lambda\in \lambda(\mu)\}$ and fix $\lambda\in \Lambda\sqcup \Lambda'$. Let $\mu_1,\mu_2 \in M_{\lambda}$. The set of étale paths between the fiber functors $\vec{\mu_1}$ and $\vec{\mu_2}$ of the category $\Rev^{\mD_{\lambda}}(\mfX_{\lambda})$ is the profinite set $\pi_1^{\mD_{\lambda}}(\mfX_{\lambda}, \vec{\mu_1},\vec{\mu_2})$ of ismorphisms between these two functors. The fundamental groups $\pi_1^{\mD_{\lambda}}(\mfX_{\lambda}, \vec{\mu_1})$ and $\pi_1^{\mD_{\lambda}}( \mfX_{\lambda}, \vec{\mu_2})$ acts by left and right translation canonically on $\pi_1^{\mD_{\lambda}}(\mfX_{\lambda}, \vec{\mu_1},\vec{\mu_2})$ and these actions are simply transitive. By construction, $\vec{\mu_1}$ and $\vec{\mu_2}$ are turned into the fiber functor $\vec{s}$ of $\Rev^{S^0}(\mS)$ through the base change by the map $f_{\lambda}\colon \mfX_{\lambda} \rightarrow \mS$ so that we have a map
		$$p_{\lambda/\mS}\colon \pi_1^{\mD_{\lambda}}(\mfX_{\lambda}, \vec{\mu_1}, \vec{\mu_2}) \longrightarrow \pi_1^{S^0}(\mS, \vec{s}).$$
		
		By composition, we get a canonical map $p_{\lambda}=p_{\mS}\circ p_{\lambda/\mS}\colon \pi_1^{\mD_{\lambda}}(\mfX_{\lambda}, \vec{\mu_1}, \vec{\mu_2})\rightarrow G_K$. 
		
		\begin{defin} 
			The geometric part $\pi_1^{\mD_{\lambda}}({\mfX_{\lambda}}_{\overline{K}}, \vec{\mu_1},\vec{\mu_2})$ of $\pi_1^{\mD_{\lambda}}(\mfX_{\lambda}, \vec{\mu_1}, \vec{\mu_2})$ is the set $p_{\lambda}^{-1}(\{1\})$. 
		\end{defin}
		
		The maps $p_{\lambda}$ (for varying $\mu\in M_{\lambda}$) induce a groupoid homomorphism from $\Pi_1^{\mD_{\lambda}}(\mfX_{\lambda}, (\vec{\mu})_{\mu\in M_{\lambda}})$ to $G_K$. This groupoid compatibility ensures that the canonical actions of the groups $\pi_1^{\mD_{\lambda}}(\mfX_{\lambda}, \vec{\mu_1})$ and $\pi_1^{\mD_{\lambda}}(\mfX_{\lambda}, \vec{\mu_2})$ on $\pi_1^{\mD_{\lambda}}(\mfX_{\lambda}, \vec{\mu_1},\vec{\mu_2})$ induce by restriction simply transitive actions from their geometric part to the geometric part of the latter. 
		
		This construction can be made when considering $\vec{\mu_1}$ and $\vec{\mu_2}$ as fiber functors with respect to the category of étale covers of $\mfX_\lambda$ tamely ramified over $\mD_\lambda'$ instead of of $\mD_{\lambda}$. As in Proposition \ref{pro:exactxlambda} we have a natural map
		\[\beta_{\mu_1,\mu_2}\colon \pi_1^{\mD_{\lambda}}(\mfX_{\lambda}, \vec{\mu_1}, \vec{\mu_2}) \longrightarrow \pi_1^{\mD_{\lambda}'}(\mfX_{\lambda}, \vec{\mu_1}, \vec{\mu_2})\]
		which is compatible with the canonical actions on both sides with regards to the maps $\beta_{\mu_1}$ and $\beta_{\mu_2}.$ In particular, the map $\beta_{\mu_1,\mu_2}$ is surjective and also induces a bijection from $\overline{p_{\lambda/\mS}}^{-1}(\{1\})$ to $\pi_1^{\mD_{\lambda}'}({\mfX_{\lambda}}_{\overline{K}}, \vec{\mu_1}, \vec{\mu_2})$ as in Proposition~\ref{pro:structfondgroup}. Moreover, the base change functor to the special fiber induces again a canonical bijection
		$$  \pi_1^{\mD_{\lambda}'}(\mfX_{\lambda}, \vec{\mu_1}, \vec{\mu_2}) \simeq \pi_1^{D_{\lambda}}(X^{0}_{\lambda}, \overrightarrow{T_{ij}}, \vec{T_{kl}}).$$
		
		\begin{defin}
			We define an action of $G_K$ on $\pi_1^{\mD_{\lambda}}(\mfX_{\lambda}, \vec{\mu_1},\vec{\mu_2})$ in the following way. For $\gamma\in \pi_1^{\mD_{\lambda}}({\mfX_{\lambda}}_{\overline{K}}, \vec{\mu_1},\vec{\mu_2})$ and $\sigma\in G_K$, let
			\[\sigma\cdot \gamma= s_{\lambda/\mu_1}(\sigma) \circ \gamma \circ s_{\lambda/\mu_2}(\sigma)^{-1}\]
			where $s_{\lambda/\mu}= p_{\mu/\lambda} \circ s_{\mu}$ for $\mu \in M_{\lambda}$. 
		\end{defin}
		
		By the compatibility with $p_{\lambda}$ this action induces an action of $G_K$ on the geometric part of $\pi_1^{\mD_{\lambda}}(\mfX_{\lambda}, \vec{\mu_1},\vec{\mu_2})$. This action is compatible with the bijection $p_{\lambda/\mS}^{-1}(\{1\})\simeq \pi_1^{D_{\lambda}}({X^{0}_{\lambda}}_{\overline{K}}, \overrightarrow{T_{ij}}, \vec{T_{kl}})$ and we recover the $G_K$-action induced by our choice of tangential base points on the right-hand side. 
		
		\subsubsection{} We can now state the groupoid analog of Proposition~\ref{pro:trivialactiongroup} and establish the main result of this section.
		
		\begin{prop}\label{pro:trivialactiongroupoid}
			Let $\mu_1,\mu_2\in M$. An element of $G_K$ acts trivially on $\pi_1^{\ell,\mD_{\lambda}}({\mfX_{\lambda}}_{\overline{K}}, \vec{\mu_1},\vec{\mu_2})$ if and only if it acts trivially on $\pi_1^{\ell,D_{\lambda}}({X^{0}_{\lambda}}_{\overline{K}}, \overrightarrow{T_{ij}}, \vec{T_{kl}}).$
		\end{prop}
		
		\begin{proof}
			As the bijection $p_{\lambda/\mS}^{-1}(\{1\})\simeq \pi_1^{D_{\lambda}}(X^{0}_{\lambda}, \overrightarrow{T_{ij}}, \vec{T_{kl}})$ is a $G_K$-isomorphism the implication is straightforward again.
			
			\medskip
			
			For the converse, let $\sigma\in G_K$. We first remark that by the simple transitiveness of the action of $\pi_1^{\mD_{\lambda}}({\mfX_{\lambda}}_{\overline{K}}, \vec{\mu_1})$ on $\pi_1^{\mD_{\lambda}}(\mfX_{\lambda}, \vec{\mu_1},\vec{\mu_2})$ and its compatibility with the map $p_{\lambda/\mS}$ we have that for every $\gamma\in \pi_1^{\ell,\mD_{\lambda}}({\mfX_{\lambda}}_{\overline{K}}, \vec{\mu_1},\vec{\mu_2})$ there exists $\alpha\in \widehat{\Z}_{\ell}(1)$ such that $\alpha\cdot \gamma \in p_{\lambda/\mS}^{-1}(\{1\})$.
			
			Now, by assumption, we have $\sigma(\alpha\cdot \gamma)= \alpha\cdot \gamma$ so that $\sigma(\gamma)=\sigma(\alpha)^{-1} \cdot (\alpha\cdot \gamma)$ and thus it is enough to see that $\sigma$ acts trivially on $\widehat{\Z}_{\ell}(1)$. This follows as in the proof of Proposition \ref{pro:trivialactiongroup}, since $\sigma$ acting trivially on $\pi_1^{D_{\lambda}}({X^{0}_{\lambda}}_{\overline{K}}, \overrightarrow{T_{ij}}, \vec{T_{kl}})$ implies it acts trivially on $\pi_1^{D_{\lambda}}({X^{0}_{\lambda}}_{\overline{K}}, \overrightarrow{T_{ij}})$, again by simple transitiveness and Galois compatibility. 
		\end{proof}
		
		\begin{rem}
			The result holds in more generality by using an almost full class of finite groups instead of the pro-$\ell$ completion.
		\end{rem}
		
		Consider the formal scheme $\mfX$. The maps $\mfX_{\lambda}\rightarrow \mfX$ for $\lambda\in \Lambda\sqcup \Lambda'$, which send $\mD$ to $\mD_{\lambda}$ by pullback, induce base change functors $\Rev^{\mD}(\mfX) \rightarrow \Rev^{\mD_{\lambda}}(\mfX_{\lambda})$. Hence for $\mu\in M$ we have fiber functors $\vec{\mu}$ for $\Rev^{\mD}(\mfX)$ and a fundamental groupoid $\Pi_1^{\mD}(\mfX, (\vec{\mu})_{\mu\in M})$ which comes with a geometric part $\Pi_1^{\mD}(\mfX_{\overline{K}}, (\vec{\mu})_{\mu\in M})$ equipped with a Galois action. For every $\lambda\in \Lambda\sqcup \Lambda'$ and $\mu_1,\mu_2\in M_{\lambda}$ the induced canonical maps
		\[p_{\lambda/\mfX,\mu_1,\mu_2}\colon \pi_1^{\mD_{\lambda}}({\mfX_{\lambda}}, \vec{\mu_1},\vec{\mu_2}) \longrightarrow \pi_1^{\mD}(\mfX, \vec{\mu_1}, \vec{\mu_2})\]
		are compatible with taking geometric parts and Galois actions on both sides.

		\begin{theo} \label{the:mainsection4}
			If an element of $G_K$ acts trivially on the groupoids $\Pi_1^{\ell,D_{\lambda}}({X^0_{\lambda}}_{\overline{K}}, \B^r_{\lambda})$ for every $\lambda\in \Lambda\sqcup \Lambda'$ then it acts trivially on the groupoid $\Pi_1^{\ell,\mD}(\mfX_{\overline{K}}, (\vec{\mu})_{\mu\in M}).$
		\end{theo} 
		
		\begin{proof}
			The main result of \cite{gm} paragraph 8.2.6 gives an equivalence of categories between $\Rev^{\mD}(\mfX)$ and a system of certain subcategories of the $\Rev^{\mD_{\lambda}}(\mfX_{\lambda})$ which yields that the fundamental groupoid  $\Pi_1^{\ell,\mD}(\mfX, (\vec{\mu})_{\mu\in M})$ is generated by the images of the $p_{\lambda/\mfX,\mu_1,\mu_2}$ for all $\lambda \in \Lambda\sqcup \Lambda'$ and $\mu_1,\mu_2 \in M_{\lambda}$. This generation statement carries to the geometric parts by \cite{iharanak} Section~3.6.
			
			The statement of the theorem now follows from Proposition \ref{pro:trivialactiongroupoid}.
		\end{proof}
		
		\medskip
		
		By Theorem~4.3.2 of \cite{gm} there is a canonical isomorphism 
		\[\Pi_1^{\ell,\mD}(\mfX, (\vec{\mu})_{\mu\in M})\simeq \Pi_1^{\ell, D}(X, (\vec{\mu})_{\mu \in M})\]
		where the right-hand side is isomorphic to $\Pi_1^{\ell}(X_{\eta}\setminus\{ (Q_v)_{v\in N}\}, (\vec{\mu})_{\mu \in M})$, and the choice of $\vec{\mu}$ defines compatible $G_K$-actions.
		
		\begin{cor}\label{cor:Xeta}
			We have the inclusion of $\ell$-monodromy fixed fields $K^{\ell}_{X_{\eta}} \subset \Q^{\ell}_{0,3}$.
		\end{cor}
		
		\begin{proof}
			For any $\mu \in M$ and $\vec{\mu}$, coming from a tangential base point of $X_{\eta}$, and seen as a fiber functor on $\Rev(X_{\eta}\setminus\{ (Q_v)_{v\in N}\}$,  we have the usual inclusion $K^{\ell}_{X_{\eta}} \subset K^{\ell}_{\vec{\mu}}$. The inclusion $K^{\ell}_{\vec{\mu}}\subset \Q^{\ell}_{0,3}$ follows by \cite{iharanak} Corollary~4.1.4~(ii). Indeed, by Theorem~\ref{the:actioncruvep1} an element of $G_K$ acts trivially on the groupoids $\Pi_1^{\ell,D_{\lambda}}({X^0_{\lambda}}_{\overline{K}}, \B^r_{\lambda})$, $\lambda\in \Lambda\sqcup\Lambda'$, if and only if it acts trivially on the groupoid $\Pi_1(\PP^1_{\overline{\Q}}\setminus \{0,1,\infty\}, \mathbb{B})$. If so, it also acts trivially on $\Pi_1^{\ell,\mD}(\mfX_{\overline{K}}, (\vec{\mu})_{\mu\in M})$ and thus on $ \Pi_1^{\ell}(X_{\eta}\setminus\{ (Q_v)_{v\in N}\}, (\vec{\mu})_{\mu \in M})$ by Theorem~\ref{the:mainsection4}.
		\end{proof}
		
		For future use let us summarize the results of Section~\ref{sec:galoisact} in a statement that can be applied for various well-chosen geometric constructions as in Section~\ref{sec:construction} of this paper. 
		
		\begin{theo}\label{theo:GM4App}
			Let $X/S$ be a stable curve with $S$ the spectrum of a discrete valuation ring with residue field $K$ of characteristic $0$. Let $D\subset X$ be a normal crossing divisor containing $X^0$ the special fiber of $X$. Let us denote by $X_{\eta}$ the generic fiber of $X$ such that $X_{\eta}$, equipped with $D_{\eta}$, is a proper smooth marked curve. Let $(X_{\lambda})_{\lambda\in \Lambda}$ be the irreducible components of $X^0$, which are equipped with a divisor $D_{\lambda}$ by pullback from $D$, and $M$ the set of double points of $X^0$. Suppose given for each $\mu\in M$ a morphism
			\[
			\mu \colon \Spf K[[T_1,T_2]]\simeq \mfX_{\mu} \to \mfX.
			\]
			If $\sigma \in G_K$ acts trivially on $\Pi_1^{D_{\lambda}}(X_{\lambda}, \{\vec{\mu}_{\lambda}\}_{\{\mu \mid \mu\ni\lambda\}})$ for every $\lambda\in \Lambda$, then it acts trivially on $\Pi_1^{D_{\eta}}(X_{\eta}, \{\vec{\mu}\}_{\mu \in M})$, where $\{\vec{\mu}_{\lambda}\}_{\{\mu \mid \mu\ni\lambda\}}$ are the associated fiber functors of $\operatorname{Rev}^{D_{\lambda}} X_{\lambda}$. 
		\end{theo}

		\section{Oda's problem for $\Z/\ell\Z$-special loci} \label{sec:conclusion}
		In the rest of this section, we fix a prime $\ell$ and specialize the previous study of this paper to the case $G=\Zl$ to establish Oda's prediction for $\Zl$-special loci -- that is the $\ell$-monodromy fixed field $\Q^{\ell}_{g,m}(\Zl)_{\underline{kr}}$ is constant independent of the topological $g$, $m$ and Hurwitz $\underline{kr}$ data and equal to $\Q^{\ell}_{0,3}$ -- which provides a new proof of Oda's original prediction, that is $\Q^{\ell}_{g,m}=\Q^{\ell}_{0,3}$.
		
		\medskip
		
		We proceed by considering two types of irreducible components $\Mgm(\Zl)_{\underline{kr}}$, whose associated monodromy fixed fields $\Q^{\ell}_{g,m}(\Zl)_{\underline{kr}}$ is compared to those of other components by the $G$-quotient of Section~\ref{sec:MonOdaG} and the $G$-deformation of Section~\ref{sec:deformation}.

		\subsection{The case of proper special loci} 
		Let us consider the case where $\MgmGkr^{\nu}$ is such that the quotient loci is $\mM_{0,3},$ that is when the quotient loci is proper. As the quotient map is itself quasi-finite and proper, the stack $\MgmGkr$ is proper if and only if it is the case of the stack of the quotient curves. In this case, both stacks $\MgmGkr$ and $\MgmGkr^{\nu}$ are geometrically given by a single point and are equal. 
		
		\medskip
		
		The following lemma enumerates the possible values of $g$, $m$ and $\underline{kr}$ that make this possible for a $\Zl$-special loci in the étale quotient case.
		
		\begin{lem}\label{lem:propCase} Assuming the ramified points are marked, the moduli space $\Mgm(\Zl)_{\underline{kr}}$ is proper in the following cases: 
			\begin{enumerate}
				\item $g=0,$ $m=2+\ell$, $\underline{k}=(1,-1)$; 
				\item $g=\frac{\ell-1}{2}$, $m=3$, and the abstract Hurwitz data $\underline{k}$ is free.
			\end{enumerate}
		\end{lem}
		
		\begin{proof}
			In the case of a quotient by $\Zl$ the Hurwitz formula is
			\[
			2g-2=(2g'-2)\ell+N (\ell-1)
			\]
			where $N$ is the number of ramified points, and setting $g'=0$ yields
			
			\[g= (N-2)(\frac{\ell-1}{2}).\]
			Since the ramified points are assumed to be marked, we have $N \in \{2,3\}$, since the cases $N=0$ or $1$ are not possible.
			
			\medskip
			
			For $N=2$ we have $g=0$ and $\underline{k}=(1,-1)$. The $m=2+\ell$ marked points are given by two ramified points and $\ell$ points permuted under the action of $\Zl$.
			
			For $N=3$ we have $g=(\ell-1)/2$ and the marked points are the ramified points. In this case, there is no condition on the abstract Hurwitz data.
		\end{proof}
		
		Let us remark that the case $N=3$ (resp. $N=2$) is given by the Seyama curves (resp. the $G$-curves of genus $0$) discussed in Section~\ref{sec:seyamacurves}. 
		
		\begin{theo}\label{the:propercase}
			For $g,m\in \N$ and compatible abstract Hurwitz data $\underline{kr}$ such that the stack $\Mgm(\Zl)_{\underline{kr}}$ is proper and non-empty, we have the equality
			\[\Q^{\ell}_{g,[m]}(\Z/\ell \Z)_{\underline{kr}}=\Q^{\ell}_{0,3}.\]
		\end{theo}
		
		Note that following our assumptions one as also $\Q^{\ell}_{g,[m]}(G)_{\underline{kr}}=\Q^{\ell}_{g,[m]}(G)_{\underline{kr}}^\nu$.

		\begin{proof}
			Corollary~\ref{cor:inclusiondiagramfields}, see diagram below, gives the inclusions $\Q^{\ell}_{0,3}\subset \Q^{\ell}_{g,[m]}(\Zl)_{\underline{kr}}\subset \Q^{\ell}_{g,[m]}(\Zl)_{\underline{kr}}^\nu$.  Let us consider $s$ and the abstract Hurwitz data $\underline{kr}^{et}$, as defined in Proposition~\ref{pro:etmap}, and the map $\MgmGkr \rightarrow \mM_{g,[m+s]}(G)_{\underline{kr}^{et}}$ which is finite. Thus $\MgmGkr$ is proper if and only if $\mM_{g,[m+s]}(G)_{\underline{kr}^{et}}$ is, and it is sufficient to establish the reverse inclusion $\Q^{\ell}_{0,3}\supset \Q^{\ell}_{g,[m]}(\Zl)_{\underline{kr}}$ in the étale quotient case, since $\Q^{\ell}_{g,[m]}(\Zl)_{\underline{kr}}\subset \Q^{\ell}_{g,[m+s]}(\Zl)_{\underline{kr}^{et}}$ by Theorem~\ref{the:nonetale}. In this case, it follows from Lemma~\ref{lem:propCase} that there is a $K$-point in the special loci that represents a curve $C$ isomorphic to either a Seyama curve or a $G$-curve of genus $0$.
			
			\medskip
			
			The result then follows from the inclusion $\Q^{\ell}_{g,[m]}(G)_{\underline{kr}}\subset\Q^{\ell}_C=\Q^{\ell}_{0,3}$ obtained from Lemma~\ref{lem:inclusioncurve} and Corollary~\ref{cor:seyamagenus0field}.
		\end{proof}
		
		\subsection{General conclusion}
		We can now establish the main result of this paper for prime cyclic special loci, which also recovers Oda's weak classical conjecture.
		
		\begin{theo}\label{theo:mainFinal}
			For $g,m\in \N$ be such that $2g-2+m>0$ and compatible abstract Hurwitz data $\underline{kr}$ such that $\Mgm(\Zl)_{\underline{kr}}$ is non-empty, we have $\Q^{\ell}_{g,[m]}(\Zl)_{\underline{kr}}= \Q^{\ell}_{0,3}$.
		\end{theo}
		
			
			\emph{Proof.} 
			By Corollary \ref{cor:nonetale} we can assume that the marked points contain the ramified points of the $G$-action. Since Theorem~\ref{the:propercase} gives the equalities $\Q^{\ell}_{g,[m]}(\Zl)_{\underline{kr}}^\nu=\Q^{\ell}_{g,[m]}(\Zl)_{\underline{kr}}= \Q^{\ell}_{0,3}$ in the case where $\Mgm(\Zl)_{\underline{kr}}$ is proper, let us assume otherwise.

			\medskip
			
				In this case, let us consider the $G$-stable diagram $X^0$ over $K$, with $\Card \Lambda\sqcup \Lambda'\geq 2$, in the boundary of $\CompMgmGkr^{\nu}$ such as provided by Proposition~\ref{pro:diagramexistence}. The stable curve $X^0$ admits a formal deformation $\mfX$ which is algebraizable into a scheme $X$ with generic fiber $X_{\eta}\in \Mgm(\Zl)_{\underline{kr}}^{\nu}(K((T)))$ as given by Section~\ref{subsub:formaldef}. The groupoid $\Pi_1^{\ell}(X_{\eta}\setminus\{ (Q_v)_{v\in N}\}, (\vec{\mu})_{\mu \in M})$ is equipped with the tangential Galois action of $G_K$ constructed in Section~\ref{sec:galoisact} coming from the choices of the fiber functors $(\vec{\mu})_{\mu \in M}$. It results from Corollary~\ref{cor:Xeta} that $K_{X_{\eta}} \subset \Q^{\ell}_{0,3}$. 
				
				\medskip
				
				It follows that $\Q^{\ell}_{g,[m]}(\Zl)_{\underline{kr}}^\nu \subset \Q^{\ell}_{0,3}$, since $\Q^{\ell}_{g,[m]}(\Zl)_{\underline{kr}}^\nu \subset K_{X_{\eta}}$ by Lemma~\ref{lem:inclusioncurve}, which concludes the first statement by the diagram below Corollary~\ref{cor:inclusiondiagramfields}. 
				In short, we obtained
				\[
				\begin{tikzcd}[column sep=10pt]
					\Q^{\ell}_{0,3} \arrow[r, hook] &
					\Q^{\ell}_{g,m} \arrow[r, hook] &
					\Q^{\ell}_{g,[m]}(\Zl)_{\underline{kr}} \arrow[r, hook] & \Q^{\ell}_{g,[m]}(\Zl)_{\underline{kr}}^\nu \arrow[r, hook] &
					K_{X_{\eta}} \arrow[r, hook] &
					\Q^{\ell}_{0,3}.
				\end{tikzcd}
				\]
				
				\qed

			Recovering Oda's weak conjecture relies on previous work of Nakamura and the consideration of certain étale type loci in $\mM_{g,[m+s]}(G)$.
			
			\begin{cor}\label{cor:newOda}
				For all $g',m'\in \N$ such that $2g'-2+m'>0$ the equality $\Q^{\ell}_{g',m'}= \Q^{\ell}_{0,3}$ holds.
			\end{cor}
			\begin{proof}
				For every $g',m'\in \N$ such that $2g'-2+m'>0$, there are $g,m\in \N$ and a compatible abstract Hurwitz data $\underline{kr}$ such that $\Mgm(\Zl)_{\underline{kr}}^{\nu}$ is non-empty and $(g',m')$ is the quotient data. This non-emptiness assertion is obtained by Proposition~3.7 of~\cite{CM15}.
				\begin{equation}\label{eq:DiagFldOda}
					\begin{tikzcd}[column sep=10pt]
						\Q^{\ell}_{g,m} \arrow[r, hook,dashed] &  \Q^{\ell}_{g,[m]}(\Zl)_{\underline{kr}} \arrow[r, hook,dashed] & \Q^{\ell}_{g,[m]}(\Zl)^{\nu}_{\underline{kr}} \arrow[r, hook]& \Q^{\ell}_{g,[m+s]}(\Zl)_{\underline{kr}^{\rm et}}^{\nu} \arrow[r, hook]& \Q^{\ell}_{0,3}\\
						\Q^{\ell}_{0,3} \arrow[u, hook, dashed] \arrow[r, hook] & \Q^{\ell}_{g',m'}\arrow[rr, hook] &  &  \Q^{\ell}_{g',m'+s'} \arrow[u, hook] & &
					\end{tikzcd}            
				\end{equation}
				From Proposition~\ref{pro:etmap} there is a non-empty stack $\mM_{g,[m+s]}(G)_{\underline{kr}^{\mathrm{et}}}$ for some $s\geq 0$ with $\underline{kr}^{\mathrm{et}}$ of étale type by construction, and such that the quotient space is $\mM_{g',m'+s'}$ for some $s'\geq 0$. By Theorem~\ref{theo:fieldQuot} we obtain the inclusion $\Q^{\ell}_{g',m'+s'}\subset \Q^{\ell}_{g,[m+s]}(G)_{\underline{kr}^{\mathrm{et}}}^{\nu}$, then $\Q^{\ell}_{g,[m+s]}(G)_{\underline{kr}^{\mathrm{et}}}^{\nu}\subset \Q^{\ell}_{0,3}$ by Theorem~\ref{theo:mainFinal}. The conclusions follows by \cite{Ue94} and \cite{Tak12} which gives the inclusion $\Q^{\ell}_{g',m'}\subset \Q^{\ell}_{g',m'+s'}$ with $s'\geq 1$, and finally by the inclusion $\Q^{\ell}_{0,3}\subset \Q^{\ell}_{g',m'}$ which is again Theorem~A of \cite{Nacoupl}, see Diag.~\ref{eq:DiagFldOda} for a summary.       
			\end{proof}
			
			\printbibliography[title={References}]
			
		\end{document}